\documentclass[11pt]{article}
  \usepackage[letterpaper, left=1in, right=1in, top=1in, bottom=1in]{geometry}
  
  \usepackage[affil-it]{authblk}  % for associating author affiliations

\def\KEYWORDSname{{\it Key words\/}{\kern0.7pt}:\enskip}
  \usepackage{enumerate}
  \usepackage{palatino}  % font face
  \usepackage{setspace}  % line spacing
  \usepackage{color}  % font color
  \usepackage[dvipsnames]{xcolor}  % adds more color options

	\usepackage[labelfont=bf]{caption}	% caption style

  \usepackage{amsmath}  % essential math package
  \usepackage{amssymb}  % math symbols
  \usepackage{amsthm}  % creating proof and theorem environments
  \usepackage{bm}  % boldface letters
  \usepackage{bbm}  % supports \mathbb-style digits, e.g. \mathbbm{1} indicating an all-one vector
  \allowdisplaybreaks[4]  % breaking multiline equations
	\usepackage{relsize}  % can increase the size of math symbols
	\usepackage{algorithm}  % floating environment for algorithmic
	\usepackage{algpseudocode}  % for easy writing of pseudo code

  \DeclareMathOperator*{\minimize}{minimize}
  \DeclareMathOperator*{\maximize}{maximize}
  \DeclareMathOperator*{\argmin}{arg\,min}
  \newcommand{\st}{\mathrm{subject\;to}}

	\theoremstyle{plain}
	\newtheorem{theorem}{Theorem}
	
	\newtheorem{lemma}{Lemma}
	\newtheorem{corollary}{Corollary}
	\theoremstyle{definition}
	
	\newtheorem{example}{Example}
	\newtheorem{assumption}{Assumption}

  \usepackage{pifont}  % dingbats and symbols

	\usepackage{longtable}  % for tables that span more than one page
	\usepackage{multirow}  % for columns spanning multiple rows
	\usepackage{array}  % for column specifications
	\usepackage{booktabs}  % providing nicer tables, using \toprule, \midrule, etc.
	\usepackage{lscape} % To output the table in landscape
  \usepackage{graphicx}  % enhanced support for graphics, allows e.g. \in­clude­graph­ics
  \usepackage{float}  % better positioning of figures, e.g. by using [h!]
  \usepackage{subfig}  % creating subfigures with \subfloat
\usepackage{pgf,tikz,pgfplots}
\pgfplotsset{compat=1.15}
\usepackage{mathrsfs}
\usepackage{epsfig}
\usepackage{epstopdf}
\usepackage{pgfplotstable}
\usepackage{ifthen}

	\usepackage{hyperref}
	\hypersetup{
		pdfstartview = {FitV},  % fits the width of the page to the window
		colorlinks = true,    % false: boxed links; true: colored links
		linkcolor = NavyBlue,  % color of internal links
		citecolor = NavyBlue,  % color of links to bibliography
		filecolor = NavyBlue,  % color of file links
		urlcolor = NavyBlue  % color of external links
	}  

  \usepackage[round]{natbib}  % for natbib, use bibtex as backend



\title{\Large{Efficient Learning of Decision-Making Models: A Penalty Block Coordinate Descent Algorithm for Data-Driven Inverse Optimization}}
\author[1]{Rishabh Gupta}
\author[1]{Qi Zhang \thanks{Corresponding author (qizh@umn.edu)}}
\affil[1]{Department of Chemical Engineering and Materials Science, \break University of Minnesota, Minneapolis, MN 55455, USA}
\date{}

\begin{document}

\maketitle

\begin{abstract}
Decision-making problems are commonly formulated as optimization problems, which are then solved to make optimal decisions. In this work, we consider the inverse problem where we use prior decision data to uncover the underlying decision-making process in the form of a mathematical optimization model. This statistical learning problem is referred to as data-driven inverse optimization. We focus on problems where the underlying decision-making process is modeled as a convex optimization problem whose parameters are unknown. We formulate the inverse optimization problem as a bilevel program and propose an efficient block coordinate descent-based algorithm to solve large problem instances. Numerical experiments on synthetic datasets demonstrate the computational advantage of our method compared to standard commercial solvers. Moreover, the real-world utility of the proposed approach is highlighted through two realistic case studies in which we consider estimating risk preferences and learning local constraint parameters of agents in a multiplayer Nash bargaining game.
\end{abstract}


\section{Introduction}

Decision making is fundamental to everyday life. As humans, we constantly make decisions that involve balancing different trade-offs to achieve the best outcome. These decisions can be as mundane as choosing what to eat for lunch, or as intricate as designing and operating a chemical plant. Also, humans are not the only ones making decisions; other biological entities, such as animals, microorganisms, and even individual cells, are often considered intelligent and autonomous agents capable of making decisions \citep{McFarland1977, Balazsi2011}. Furthermore, with increased automation and the emergence of artificial intelligence, many engineered systems can also be viewed as decision-making agents \citep{Steels1995}.

A good understanding of decision-making mechanisms is crucial for predicting the behavior of autonomous agents, learning from experts, and optimizing systems involving various decision makers. But many decision-making processes are unknown or poorly understood. For example, experts in the operation of chemical plants make decisions based on many years of experience, but their decision strategies often are not well documented and, due to the complexity of the manufacturing processes, are difficult to explain even to fellow operators. As a result, the complete transfer of expert knowledge to new operators remains an unsolved problem. Likewise in microbiology, cells can be considered autonomous agents that make decisions regarding gene expression and cell metabolic function. While we can observe these decisions in experiments, we often do not understand why cells make these choices. Answering this question would provide fundamental insights that could advance cancer treatment, immunology research, and biomanufacturing operations. In both examples, we do not know how exactly decisions are made, but we can observe those decisions. The question is: Can we use these observations to uncover the underlying decision-making process?

The above question is, of course, not new. In particular, the first motivating example may remind some readers of expert systems, a subdomain of artificial intelligence that had spurred much enthusiasm in various disciplines, including chemical engineering, in the 1980s \citep{Banares-Alcantara1985, Rich1987, Stephanopoulos1990}. While research in this field has primarily focused on using human experience and domain knowledge to build expert systems, there have been a few attempts to learn an expert system from data \citep{Rich1989, Sammut1992}, most of which are based on decision tree learning. More recently, in machine learning, this problem has been referred to as apprenticeship or imitation learning \citep{Abbeel2004}, with popular applications in robotics and autonomous driving. In theory, any machine learning method can be applied to this problem, but it is unclear what model is the most suitable. Popular black-box machine learning models, such as artificial neural networks, often require a large amount of data to train and, more importantly, are difficult to interpret. Interpretability, however, is crucial in this problem as we are chiefly interested in gaining a better understanding of the underlying decision-making process.

In this work, we present an approach that is fundamentally inspired by the \textit{principle of optimality} \citep{Schoemaker1991}, which conjectures that autonomous agents generally make decisions in some optimal fashion. This basic relationship between decision making and optimization is commonly applied in many fields. For example, most economic theory relies on the assumption that rational decision makers are utility maximizers \citep{Morgenstern1944, Hey1994}. Evolutionary biology tells us that through adaptations over the course of millions of years, biological systems have evolved to behave optimally, albeit the measure of optimality and the decision set often being unclear \citep{Rosen1967, Parker1990}. Even the fundamental basis for thermodynamics is an optimization problem, namely Gibbs free energy minimization \citep{Gautam1979, Rossi2009}. Following this principle of optimality, the key idea is to model a decision-making process as a mathematical optimization problem of the following general form:
\begin{equation*}
\begin{aligned}
\minimize_x \quad & f(x,u) \\
\st \quad & g(x,u) \leq 0,
\end{aligned}
\end{equation*}
where $x$ and $u$ denote the vectors of decision variables and input parameters, respectively. Given input parameters $u$, the decision maker chooses $x$ such that their objective function $f(x,u)$ is minimized while satisfying the constraints $g(x,u) \leq 0$. We then consider the inverse problem, which is to infer the functions $f$ and $g$ given observations, with each data point $i$ corresponding to an input-decision pair $(u_i,x_i)$. The resulting estimated model will be a natural representation of a decision-making process that is inherently interpretable. Moreover, this approach allows us to take advantage of all the modeling flexibility provided by mathematical optimization and directly incorporate domain knowledge in the form of constraints. In the operations research literature, the inverse problem described above is referred to as \textit{inverse optimization} (IO) \citep{Ahuja2001}. We present a review of the relevant IO literature in the next section.

To better understand when IO will produce a more useful model than a black-box machine learning method, consider a decision maker who is responsible for routing traffic on a road network. Our goal is to understand the decision-making strategy of this decision maker because they are considered an expert. When using IO, the constraints on the decision maker, which represent the structure of the network and its capacities, can be easily formulated. Through our approach, learning the objective function of the decision maker will uniquely allow us to understand the preferences of the expert decision maker. This information can be used to transfer the expert's knowledge to other novice operators, which can be particularly useful in situations where the expert is unavailable or when their expertise needs to be shared with a larger group of decision makers. It is generally not easy to extract such insights from black-box machine learning models.

\section{Literature Review}
\label{sec:LitReview}

The notion of inverse optimization was first introduced by \citet{Burton1992} who considered the problem of recovering arc costs on a directed graph using given shortest-path solutions. Early works have primarily addressed the deterministic setting in which observations are assumed to be exactly optimal solutions of the optimization model, with the majority of existing works only considering one single observation \citep{Zhang1996b, Zhang1996, Ahuja2001, Heuberger2004, Iyengar2005}. 

IO received relatively little attention in the 1990s and 2000s. However, there has been renewed interest in recent years when the research focus has shifted towards the case with multiple observations, also referred to as \textit{data-driven} IO \citep{Esfahani2018}, and noisy data \citep{Chan2014, Chan2019, Keshavarz2011, Bertsimas2015, Aswani2018}. This new perspective has significantly broadened the settings to which IO can be applied and it has found application in a diverse range of fields, including healthcare planning \citep{Chan2014, Chan2021a}, transportation network design \citep{Chow2014, Ronnqvist2017}, electricity markets \citep{Saez-Gallego2016, Birge2017}, biological network inference \citep{Burgard2003, uygun2007, Zhao2015}, and expert systems \citep{Akhtar2022}; a more detailed description of existing IO applications can be found in \cite{Chan2021}.

While initially designed to determine an objective function that renders a single given solution optimal, the recent research in data-driven IO has helped it gain increased acceptance as a statistical learning method \citep{Iraj2021}. Nevertheless, it has received virtually no attention from the chemical engineering or, more specifically, the process systems engineering (PSE) community. We could only find a handful of articles where IO (or similar) methods have been used to infer parameters of an optimization problem in PSE \citep{Burgard2003, uygun2007, Bollas2009, Glass2017a, Glass2018}. We believe that IO can potentially be very effective in PSE applications as it provides a compelling framework for the integration of optimization, machine learning, and domain knowledge. This is in the same spirit as other hybrid modeling methods that incorporate first-principles knowledge into otherwise purely data-driven approaches, which is crucial in many chemical engineering applications \citep{Boukouvala2017, Wilson2019, Bangi2020}.

From a theoretical perspective, research on data-driven IO is largely limited to the case where the decision-making process is being modeled as a convex optimization problem. Here, the main distinction among the various proposed formulations is in terms of the loss function employed to fit the data. Minimization of the slack required to make the noisy data satisfy an optimality condition is considered in \cite{Boyd2011}, \cite{Bertsimas2015}, and \cite{Esfahani2018}. \cite{Aswani2018} showed that this kind of loss function can lead to statistically inconsistent estimates and proposed minimizing the sum of some norm of residuals with respect to the decision variables. 

A key challenge in using the statistically consistent bilevel IO problem formulation proposed by \cite{Aswani2018} is its computational intractability. This has restricted the use of IO to very specific problems for which efficient solution methods have been developed. Almost all existing literature addresses inverse linear optimization problems \citep{Babier2021, Shahmoradi2021a, gupta2022deco}. Additionally, most of these studies assume the constraints to be known and only attempt to estimate the cost coefficients. There are very few existing contributions that consider joint estimation of objective and constraint parameters \citep{Chan2019a, Ghobadi2020}. Hence, there is a need for more efficient algorithms that can infer all unknown parameters for high-dimensional convex nonlinear problems with large datasets. This work aims to fill this apparent deficiency in the literature.

% {\color{Maroon} (I feel we don't need this last part.) Most of these studies assume the parametric form of the cost function to be known and attempt to find a set of parameters that can explain the observed behavior most closely. \cite{Bertsimas2015} is the only study that considers non-parametric estimation of the cost function.} 

\section{Outline and Contributions}
In this work, we consider data-driven IO with noisy data, and we specifically focus on the case where the optimization model describing the decision-making process can be assumed to be convex. Here, one major challenge is the computational complexity of the IO problem as it gives rise to a bilevel optimization problem whose size increases with the number of data points. To solve large model instances, we develop a penalty block coordinate descent (BCD) algorithm that exploits the specific decomposable structure of the problem. The efficacy of the proposed algorithm is demonstrated through a comprehensive set of computational experiments. Moreover, importantly, the computational results show empirically that the proposed data-driven IO method is statistically consistent and data-efficient, highlighting its promise for future, more complex applications. 

In summary, the key contributions of this work are as follows:
\begin{enumerate}
\item We propose a data-driven framework to simultaneously learn the objective function and constraint parameters of a convex nonlinear optimization problem based on multiple noisy observations of its optimal solutions for different input parameter values.
\item We formulate the data-driven IO problem as a bilevel program, which we further reformulate into a single-level optimization problem using the Karush-Kuhn-Tucker (KKT) optimality conditions. In general, the resulting problem lacks regularization. We alleviate this deficiency by proposing an exact penalty reformulation.
\item For certain important classes of functions $f$ and $g$, we show that the penalty reformulation of the IO problem is a multiconvex optimization problem. In this case, we show that large instances of the reformulated IO problem can be solved very efficiently by exploiting their decomposable structure with the BCD algorithm.
\item By applying known results from the literature, we show that the proposed solution method is guaranteed to converge to a stationary point of the reformulated IO problem. We further demonstrate the effectiveness of our method through computational experiments based on a number of instances of a convex optimization problem from the wireless communication literature. Our results indicate that a BCD-based decomposition strategy is significantly faster and finds higher-quality solutions than standard nonlinear optimization solvers.
\item We demonstrate the real-world utility of IO by considering two realistic example problems. In the first problem, we use IO to estimate a decision maker's risk preference when making decisions under significant uncertainty. In the second example, we use IO to estimate customers' internal constraints under the assumption that their buying decisions are covered by a collective bargaining agreement.
\end{enumerate}
The remainder of this paper is organized as follows. In Section \ref{sec:ProblemStatement}, we present the mathematical formulation of the IO problem that we consider in this work. The penalty BCD algorithm is developed in Section \ref{sec:Algorithm}, and the results from the computational case studies are presented in Section \ref{sec:CompRes}. Finally, we close with some concluding remarks in Section \ref{sec:Conclusions}.

\section{The Inverse Optimization Problem}
\label{sec:ProblemStatement}

In this section, we present the mathematical formulation of the IO problem that we consider in this work. We describe its key properties and introduce a set of assumptions that will be useful in the next section where we present our solution strategy for the problem.

We start by assuming that the given decision-making process can be modeled as a convex optimization problem of the following form:

\begin{equation}
\label{eqn:FOP}
\tag{FOP}
\begin{aligned}
  \minimize_{x \in \mathbb{R}^n} \quad & f(x,u;\theta) \\
  \st \quad & g_k(x,u;\omega) \leq 0 \quad \forall \, k \in \{ 1, \ldots, m \} \\
  & h_k(x, u; \omega) = 0 \quad \forall \, k \in \{ 1, \ldots, p \},
\end{aligned}
\end{equation}
where $x$ is the vector of decision variables and $u$ is the vector of input variables. The function $f$ is the objective function of the problem, and $g_k$ and $h_k$ are the constraint functions. The model given by $f$, $g_k$, and $h_k$, which we refer to as the forward optimization problem (FOP), is parameterized by two sets of model parameters, $\theta$ and $\omega$. To ease the notation, we denote the constraint functions $\{g_k\}_{k = 1}^m$ and $\{h_k\}_{k = 1}^p$ by function vectors $g$ and $h$, respectively, in the remainder of this paper.

To ensure convexity of \eqref{eqn:FOP}, we require the functions $g$ and $h$ to be convex and affine in $x$, respectively, for fixed $u$ and $\omega$. We also make the following assumptions.

\begin{assumption} \label{asp:diff}
The functions $f(x, u; \theta)$ and $g(x, u; \omega)$ are twice continuously differentiable in $x$ for fixed $u$, $\theta$, and $\omega$.
\end{assumption}

\begin{assumption} \label{asp:stconv}
The function $f(x,u;\theta)$ is \textit{strictly convex} in $x$ for fixed $u$ and $\theta$.
\end{assumption}

Assumption \ref{asp:diff} is fairly standard in the IO literature \citep{Boyd2011, Bertsimas2015, MohajerinEsfahani2018} and is required to ensure that the bilevel IO problem (introduced later) can be reformulated as a single-level problem. There are two main arguments for Assumption 2:
\begin{enumerate}
    \item Given a dataset consisting of $(u_i, x_i)$ observations, it is likely that there will be multiple optimization models that can explain the data. In these cases, it may be better to learn a strictly convex model that, when used for predictions, results in unique decisions for every input. This can help avoid ambiguities and ensure that the model makes clear, distinct decisions.
    \item There may still be situations in which the user wants to learn parameters of a problem that, generally, does not have unique solutions. However, as shown in our previous work \citep{gupta2022deco}, methods developed for estimating parameters of generally convex FOPs do not provide good solutions when the FOPs can admit multiple solutions. Therefore, in this case, we suggest adding a small regularizing quadratic term to the objective function of the problem to make it strictly convex. 
    \end{enumerate}

% \end{itemize}}Assumption \ref{asp:stconv} is motivated by our  

The goal of IO is to recover model \eqref{eqn:FOP} from observations that are assumed to be noisy estimates of the optimal solutions to \eqref{eqn:FOP} for given values of the input variables. Specifically, the problem we aim to solve is the following \citep{Aswani2018}:
\begin{equation}
\label{eqn:IOP}
\tag{IOP}
\begin{aligned}
  \minimize_{\hat{\theta} \in \Theta, \, \hat{\omega} \in \Omega, \, \hat{x}} \quad & \sum_{i \in \mathcal{I}} (x_{i} - \hat{x}_{i})^{\top} W (x_{i} - \hat{x}_{i}) \\
  \st \quad & \hat{x}_{i} = \arg \min_{\tilde{x} \in \mathbb{R}^n} \left\lbrace f(\tilde{x},u_i;\hat{\theta}): g(\tilde{x},u_i;\hat{\omega}) \leq 0, \, h(\tilde{x}, u_i; \hat{\omega}) = 0  \right\rbrace \quad \forall \, i \in \mathcal{I},
\end{aligned}
\end{equation}
where $\mathcal{I}$ is the set of experiments, with each experiment $i$ defined by given inputs $u_i$ and the corresponding observed decisions $x_i$. The objective is to choose the parameters $\hat{\theta}$ and $\hat{\omega}$ from sets $\Theta$ and $\Omega$, respectively, such that a weighted sum of squared residuals, where $W \in \mathbb{R}^{n \times n}$ denotes the diagonal matrix of weighting factors, is minimized. 

We make the following assumption about sets $\Theta$ and $\Omega$, which is required to ensure statistical consistency of the estimates $\hat{\theta}$ and $\hat{\omega}$ obtained from \eqref{eqn:IOP} \citep{Aswani2018}:

\begin{assumption} \label{asp:comp}
The sets $\Theta$ and $\Omega$ are compact.
\end{assumption}

In addition, we require a regularity assumption on these sets to ensure the convergence of our solution algorithm for \eqref{eqn:IOP}.

\begin{assumption} \label{asp:regular_sets}
The description of sets $\Theta$ and $\Omega$ satisfies necessary regularity conditions such as Mangasarian-Fromovitz constraint qualification (MFCQ).
\end{assumption}

Finally, we highlight an implicit assumption while formulating an IO problem where we assume that the set of possible objective parameters $\Theta$ disallows the trivial solution for $\theta$ that would render the objective function a constant.
% Convexity of the sets is required for the convergence of the solution algorithm proposed in this paper. 

Problem \eqref{eqn:IOP} is similar to how the IO problem has typically been formulated in the literature with the only difference being the weighting factors that offer the following flexibilities while formulating the IO problem:
\begin{enumerate}
    \item Weights can be used to reflect the level of importance of each variable. In data reconciliation, the weight for each variable is the reciprocal of the variance of the corresponding measurement, which reflects the accuracy of the measurement. This can be considered using our general matrix of weighting factors.
    \item Weighting factors allow us to deal with problems where different decision variables have vastly different scales.
    \item It also generalizes the case with unmeasured variables, for which the weighting factors can be set to zero.
\end{enumerate}
 
Finally, we close this section by formally stating a mild assumption under which \eqref{eqn:IOP} is always feasible.

 \begin{assumption} \label{asp:FeasFOP}
The sets $\mathcal{C}_i := \left\{x: g(x, u_i; \omega) < 0, h(x, u_i; \omega) = 0 \right\} \; \forall \, i \in \mathcal{I}$ are nonempty for every $\omega \in \Omega$, and $f$ remains bounded from below for every $\theta \in \Theta$.
 \end{assumption}

%  \begin{proposition} \label{prop:FeasIOP}
%  If Assumption \ref{asp:FeasFOP} holds, \eqref{eqn:IOP} is always feasible.
%  \end{proposition}
%  \begin{proof}
%  Proposition \ref{prop:FeasIOP} follows trivially from Assumption \ref{asp:FeasFOP}.
%  Let $x^*(u;\theta, \omega)$ be the optimal solution of \eqref{eqn:FOP} for input parameters $u$ and model parameters $\theta$ and $\omega$.  Without loss of generality, assume that the sets $\Theta$ and $\Omega$ contain $N$ different values of the parameters $\theta$ and $\omega$ respectively; let these values be denoted by $\{\theta_p, \omega_p\}_{p = 1}^N.$ In this case, \eqref{eqn:IOP} can be solved by first obtaining $x^*(u_i; \theta_p, \omega_p)$ for all $i \in \mathcal{I}$ and $p = \{1, \ldots, N\},$ and later calculating the objective value for each of the $N$ different model parameter combinations with $\hat{x}_i = $
%  \end{proof}

\section{Solution Strategy}
\label{sec:Algorithm}

\subsection{Single-Level Reformulation}

Problem \eqref{eqn:IOP} is a bilevel optimization problem with convex lower-level problems. This class of problems are typically solved by reformulating them into a single-level optimization problem by replacing the lower-level problems with their optimality conditions \citep{dempe2015bilevel}. Here, given Assumption \ref{asp:FeasFOP}, the lower-level problems of \eqref{eqn:IOP} can be replaced with their KKT conditions, resulting in the following reformulation of \eqref{eqn:IOP}:
% \begin{assumption} \label{asp:Slater}
% The sets $\mathcal{C}_i := \left\{x: g(x, u_i; \omega) < 0, h(x, u_i; \omega) = 0 \right\} \; \forall \, i \in \mathcal{I}$ are nonempty.
% \end{assumption}
% Note that this assumption is mild as the feasible set for some $u_i$ that does not satisfy this property can always be perturbed slightly to ensure that the assumption gets satisfied. 
\begin{subequations}
\label{eqn:KKT-Convex}
    \begin{align}
    \minimize\limits_{\hat{\theta} \in \Theta, \, \hat{\omega} \in \Omega, \hat{x} ,\lambda, \mu} \quad & \sum\limits_{i \in \mathcal{I}}(x_{i} - \hat{x}_{i})^{\top} W (x_{i} - \hat{x}_{i}) \\
    \st \quad \; \; &  \nabla {f} (\hat{x}_i, u_i; \hat{\theta}) + \lambda_i^{\top} \nabla {g} (\hat{x}_i, u_i; \hat{\omega}) + \mu_i^{\top} \nabla {h} (\hat{x}_i, u_i; \hat{\omega}) = 0 \quad \forall \, i \in \mathcal{I} \label{eqn:stat} \\
    & {g}(\hat{x}_i, u_i; \hat{\omega}) \leq 0 \quad \forall \, i \in \mathcal{I} \label{eqn:PF1} \\
    & {h}(\hat{x}_i, u_i; \hat{\omega}) = 0 \quad \forall \, i \in \mathcal{I} \label{eqn:PF2} \\
    & \lambda_{ik} \, {g}_k(\hat{x}_i, u_i; \hat{\omega}) = 0 \quad \forall \, i \in \mathcal{I}, k \in \mathcal{K}  \label{eqn:CS} \\
    & \lambda_i \geq 0, \hat{x}_i \in \mathbb{R}^n \quad \forall \, i \in \mathcal{I},
    \end{align}
\end{subequations}
where $\lambda$ and $\mu$ are respectively the dual variables of the inequality and equality constraints of the lower-level problems of \eqref{eqn:IOP}. Constraints \eqref{eqn:stat}, \eqref{eqn:PF1}-\eqref{eqn:PF2}, and \eqref{eqn:CS} correspond to the stationarity, primal feasibility, and complementary slackness conditions, respectively. 

Problem \eqref{eqn:KKT-Convex} is generally a nonconvex nonlinear program (NLP) that is difficult to solve to global optimality. Moreover, the size of \eqref{eqn:KKT-Convex} increases with the number of data points, which quickly leads to very large instances in problems of practical relevance. In the machine learning literature, given a large nonconvex optimization problem, it is typical to focus on finding a good local solution efficiently rather than trying to solve the problem to global optimality. A popular example is deep learning, which is a nonconvex optimization problem; however, several local solution methods have been found to result in very accurate predictive models \citep{choromanska2015loss, Jin2021, danilova2022}. In the following, we pursue a similar strategy in tackling the computational complexity of problem \eqref{eqn:KKT-Convex}.

\subsection{Penalty Reformulation} 

Problem \eqref{eqn:KKT-Convex} lacks sufficient regularity that is generally required to solve such problems with standard NLP methods. It contains complementarity constraints \eqref{eqn:CS}, which make the problem violate the MFCQ everywhere in the feasible region \citep{scheel2000mathematical, anitescu2005using}. Moreover, the unknown FOP parameters $\theta$ and $\omega$ complicate the form of the other constraints in \eqref{eqn:KKT-Convex} that may, in some cases, also contribute towards the lack of regularization. 

The aforementioned lack of regularization of \eqref{eqn:KKT-Convex} is known to result in degenerate and unbounded Lagrange multipliers in the vicinity of the optimal solutions, causing convergence difficulties with NLP solution algorithms. A popular strategy to overcome this difficulty employs a penalty reformulation \citep{bertsekas1997nonlinear, nocedal2006numerical} of the original problem \citep{anitescu2000nonlinear}. In this approach, the complicating constraints, i.e., the constraints resulting in constraint qualification violation are moved to the objective function as penalty terms, leaving the problem with a ``regular" feasible region. Depending on the type of penalty function, the reformulation can be exact in the sense that, for certain finite values of the penalty parameters, every optimal solution of the original problem will also be optimal for the reformulation. Commonly employed examples of penalty functions that result in these exact reformulations are nonsmooth penalty functions, such as the one based on the $\ell_1$-norm of constraint violation \citep{bertsekas1997nonlinear, nocedal2006numerical}. These have proven to be very successful in dealing with the challenges associated with difficult NLPs such as mathematical programs with complimentarity constraints (MPCCs) \citep{benson2006interior}, which motivates us to consider the following reformulation of \eqref{eqn:KKT-Convex}:   

\begin{equation}
\label{eqn:KKT-EPR}
    \begin{aligned}
    \minimize\limits_{\hat{\theta} \in \Theta, \, \hat{\omega} \in \Omega, \hat{x} ,\lambda, \mu} \quad & \sum\limits_{i \in \mathcal{I}} (x_{i} - \hat{x}_{i})^{\top} W (x_{i} - \hat{x}_{i}) 
    \\ & + c^{\top} 
    \underbrace{\begin{bmatrix} 
    \sum\limits_{i \in \mathcal{I}} \lvert \nabla \hat{f} (\hat{x}_i, u_i; \hat{\theta}) + \lambda_i^{\top} \nabla \hat{g} (\hat{x}_i, u_i; \hat{\omega}) + \mu_i^{\top} \nabla \hat{h} (\hat{x}_i, u_i; \hat{\omega}) \rvert \\
    \sum\limits_{i \in \mathcal{I}} \mathrm{max} \{ 0, \hat{g}(\hat{x}_i, u_i; \hat{\omega}) \} \\
    \sum\limits_{i \in \mathcal{I}} \lvert \hat{h}(\hat{x}_i, u_i; \hat{\omega}) \rvert \\
    \sum\limits_{i \in \mathcal{I}} \lvert \lambda_{i}^{\top} \, \hat{g}(\hat{x}_i, u_i; \hat{\omega}) \rvert
    \end{bmatrix}}_{P} \\
    \st \quad & \lambda_i \geq 0, \hat{x}_i \in \mathbb{R}^n \quad \forall \, i \in \mathcal{I} \\
    \end{aligned}
\end{equation}
where $c$ are the positive penalty parameters. 

Next, we present our strategy to find a minimizer for \eqref{eqn:KKT-Convex} through \eqref{eqn:KKT-EPR}.

\begin{lemma} \label{lem:regular}
The feasible region of problem \eqref{eqn:KKT-EPR} satisfies MFCQ.
\end{lemma}
\begin{proof}
By Assumption \ref{asp:regular_sets}, the constraints describing the sets $\Theta$ and $\Omega$ satisfy MFCQ. Moreover, in \eqref{eqn:KKT-EPR}, the only remaining constraints are the non-negativity constraints on the dual variables $\lambda$ which satisfy linear independence constraint qualification (LICQ) and therefore MFCQ.
\end{proof}

In \eqref{eqn:KKT-EPR}, we move all but a few simple bound constraints of \eqref{eqn:KKT-Convex} to the objective function. However, note that our goal here is to only eliminate the constraints that contribute to the lack of regularity of \eqref{eqn:KKT-Convex}. Therefore, while formulating \eqref{eqn:KKT-EPR}, one can choose to move only those constraints for which constraint qualification violation is expected. Since, in practice, it is often not trivial to check constraint qualification violation, we use the following reformulation strategy: We always penalize all the stationarity, complementarity, and primal feasibility constraints for which $\omega$ are unknown. For the other primal constraints with known $\omega$, we leave them as hard constraints if it is possible to see that they do not violate constraint qualification. For example, if the remaining primal constraints are simple bound constraints that obviously satisfy LICQ, then we do not need to penalize them.

% The following theorem \citep{nocedal2006numerical} relates the solutions of \eqref{eqn:KKT-Convex} with those of \eqref{eqn:KKT-EPR}:

% \begin{theorem} \label{thm:Exactness}
% Suppose $(\hat{\theta}^*, \hat{\omega}^*, \hat{x}^*, \lambda^*, \mu^*)$ minimizes \eqref{eqn:KKT-EPR} for any $c \geq \tilde{c}$. If $(\hat{\theta}^*, \hat{\omega}^*, \hat{x}^*, \lambda^*, \mu^*)$ is feasible to \eqref{eqn:KKT-Convex}, then it also minimizes \eqref{eqn:KKT-Convex}.
% \end{theorem}

Lemma \ref{lem:regular} essentially highlights the well-posedness of \eqref{eqn:KKT-EPR} which makes it amenable to solution via any commercial NLP solver, such as IPOPT \citep{wachter2006implementation}. We exploit this property of \eqref{eqn:KKT-EPR} to design an algorithm in which the solution of \eqref{eqn:KKT-Convex} can be obtained by repeatedly solving \eqref{eqn:KKT-EPR}. The strategy employed in Algorithm \ref{alg:Outerloop} is to start with some initial values for the penalty parameters $c$ and successively increase them by a factor $\rho$ after every iteration to find the point where a feasible solution to \eqref{eqn:KKT-Convex} is obtained. The global convergence of this approach for a general NLP has been proved by \cite{benson2009convergence}.

\begin{algorithm}
\begin{algorithmic}[1]
    \State initialize: $ k \gets 1, (\hat{\theta}, \hat{\omega}, \hat{x}, \lambda, \mu) \gets (\hat{\theta}_0, \hat{\omega}_0, \hat{x}_0, \lambda_0, \mu_0)$ and $c \gets c_1$
    \While{$\lVert P \rVert > \epsilon $}
        \State solve \eqref{eqn:KKT-EPR} with $c = c_k$, warm-start with $(\hat{\theta}_{k-1}, \hat{\omega}_{k-1}, \hat{x}_{k-1}, \lambda_{k-1}, \mu_{k-1})$, obtain $(\hat{\theta}_{k}, \hat{\omega}_{k}, \hat{x}_{k}, \lambda_{k}, \mu_{k})$
        \State $c_{k+1} \gets c_k + \rho c_k$
        \State $k \gets k+1$
    \EndWhile
    \State \Return $(\hat{\theta}_k, \hat{\omega}_k, \hat{x}_k, \lambda_k, \mu_k)$
    \caption{Algorithm for solving \eqref{eqn:KKT-Convex} through \eqref{eqn:KKT-EPR}.}
    \label{alg:Outerloop}
\end{algorithmic}
\end{algorithm}

\begin{assumption} \label{asp:bounded}
Suppose $(\hat{\theta}_k, \hat{\omega}_k, \hat{x}_k, \lambda_k, \mu_k)$ denotes a solution of \eqref{eqn:KKT-EPR} with $c = c_k$. The sequence of solutions $\{\hat{\theta}_k, \hat{\omega}_k, \hat{x}_k, \lambda_k, \mu_k\}$ generated as $c_k$ increase to $\infty$ is bounded.
\end{assumption}

\begin{theorem} \citep{benson2009convergence}
Suppose Assumption \ref{asp:bounded} holds. If $(\hat{\theta}^*, \hat{\omega}^*, \hat{x}^*, \lambda^*, \mu^*)$ is the accumulation point of the sequence of solutions generated as $k \rightarrow \infty$ and $c_k$ increase to $\infty$, then $(\hat{\theta}^*, \hat{\omega}^*, \hat{x}^*, \lambda^*, \mu^*)$ minimizes \eqref{eqn:KKT-Convex}.
\end{theorem}

% Finally, it is worth noting here that the existence of a finite $\bar{c}$ for which \eqref{eqn:KKT-EPR} is an exact reformulation of \eqref{eqn:KKT-Convex} typically requires \eqref{eqn:KKT-Convex} to satisfy certain constraint qualification such as MFCQ as shown in \citep{han1979exact}; this is clearly not the case here as highlighted before. However, it has been shown in several cases, such as when \eqref{eqn:KKT-Convex} satisfies a relaxed CQ known as MPEC-LICQ

\subsection{Block Coordinate Descent}
We introduced the penalty reformulation as a means to tackle the lack of regularity of \eqref{eqn:KKT-EPR}. However, the resulting ``well-posed" problem is still a nonconvex NLP whose size increases with the number of data points. In this section, we propose a decomposition scheme designed to significantly reduce the computation time when solving certain large instances of \eqref{eqn:KKT-EPR}.

Our approach is motivated by the fact that for many FOPs, \eqref{eqn:KKT-Convex} becomes a \textit{multiconvex optimization problem} (MCP) that exhibits a particular decomposable structure. MCPs are problems for which the variables can be partitioned into blocks over which it is convex when all other variables are held constant \citep{shen2017disciplined}. Such FOPs include, for example, the broad class of quadratic programs (QPs), which are very common in engineering applications \citep{mccarl1977}. Next, we demonstrate the MCP structure of the penalty reformulation by deriving \eqref{eqn:KKT-EPR} for the case when \eqref{eqn:FOP} is a QP.

\begin{example}
Suppose that \eqref{eqn:FOP} can be posed as the following convex QP:
\begin{equation}
    \begin{aligned}
    \minimize\limits_{x \in \mathbb{R}^n} \quad & \frac{1}{2} x^{\top}Qx + r^{\top}x \\
    \st \quad & Ax \leq b,
    \end{aligned}
\end{equation}
where $Q \succ 0$. We assume that $Q$ and $r$ can be parameterized with $\theta$, and the parameters $A$ and $b$ can be parameterized with $\omega$. Also, $Q$, $r$, $A$, and $b$ generally change with the input parameters $u$. Problem \eqref{eqn:KKT-EPR} for this FOP can be stated as follows:
\begin{equation}
\begin{aligned}
    \minimize_{\hat{\theta} \in \Theta, \, \hat{\omega} \in \Omega, \hat{x} , \lambda} \quad & \sum_{i \in \mathcal{I}} (x_{i} - \hat{x}_{i})^{\top} W (x_{i} - \hat{x}_{i}) \\
&+ c^{\top} 
\begin{bmatrix} 
  \displaystyle  \sum_{i \in \mathcal{I}} \lvert Q (u_i;\hat{\theta}) \hat{x}_i + r (u_i;\hat{\theta}) + \lambda^{\top}_i A (u_i;\hat{\omega}) \rvert \\
  \displaystyle  \sum_{i \in \mathcal{I}} \lvert A(u_i;\hat{\omega}) \, \hat{x}_i - b (u_i;\hat{\omega}) \rvert \\
    \displaystyle\sum_{i \in \mathcal{I}} \lvert \lambda_{i}^{\top} \left( A(u_i;\hat{\omega}) \, \hat{x}_i - b (u_i;\hat{\omega}) \right) \rvert
\end{bmatrix}
 \\
    \st \quad & \lambda_{i} \geq 0, \hat{x}_i \in \mathbb{R}^n \quad \forall \, i \in \mathcal{I}.
\end{aligned}
\end{equation}
Assuming that $\Theta$ and $\Omega$ are convex, this problem is an MCP for the following block decomposition of its variables: $\left((\hat{\theta}, \lambda), \hat{\omega}, \hat{x}\right).$
\end{example}

MCPs are nonconvex problems and hence generally very hard to solve globally; however, several important problems in machine learning fall into this class and have led to the development of a number of approximate solution methods \citep{shen2017disciplined, jain2017non}. Most of these methods employ some variation of the \textit{block coordinate descent} (BCD) approach which exploits the particular mathematical structure of subproblems obtained by fixing blocks of variables. 

Assume that the MCP instances of \eqref{eqn:KKT-EPR} can be written as follows:
\begin{equation}
    \label{eqn:MCP}
    \begin{aligned}
    \minimize\limits_{x} \quad & f(x_1, \ldots, x_\ell) \\
    \st \quad & x_i \in \mathcal{X}_i \quad \forall \, i \in \{1, \ldots, \ell \},
    \end{aligned}
\end{equation}
where the variables $x \in \mathbb{R}^n$ have been partitioned into $\ell$ mutually exclusive, collectively exhaustive sets $\{x_1, \ldots, x_\ell\}.$ The function $f$ is the nonconvex objective function and $\mathcal{X}_i$ is a closed convex set of feasible points for the variable block $x_i$. 
Note that in \eqref{eqn:MCP}, we require different variable blocks to lie in their own independent feasible sets; this is a requirement for the BCD method as it has been shown that in the presence of coupled feasible sets, the BCD algorithm may stagnate at a non-stationary point. Here, our penalty reformulation naturally gives us a formulation \eqref{eqn:KKT-EPR} for which the feasible sets for different variables are independent of each other. 

Problem \eqref{eqn:MCP} is an MCP because the following subproblems are convex:
\begin{equation}
    \label{eqn:SubMCP}
    \begin{aligned}
    \minimize\limits_{{x}_i \in \mathcal{X}_{i}} \quad & f(\bar{x}_1, \ldots, \bar{x}_{i-1}, x_{i}, \bar{x}_{i+1}, \ldots, \bar{x}_\ell) \quad \forall \, i \in \{1, \ldots, \ell\}, \\
    \end{aligned}
\end{equation}
where the bar over a variable indicates a fixed variable. In BCD, one solves \eqref{eqn:MCP} by solving the subproblems \eqref{eqn:SubMCP} cyclically in a Gauss-Seidel fashion until the variable values converge to a limit point. While the required number of iterations for BCD to converge can be exponentially large, the problems solved at each step are convex and overall, this approach has been proven to be significantly more efficient than directly solving the full-space problem for many practical problems \citep{Wright2015}. 

Generally, the convergence of BCD to a stationary point of an optimization problem requires rather stringent conditions on the structure of \eqref{eqn:MCP}; typically, $f$ needs to be smooth and all the subproblems \eqref{eqn:SubMCP} must have unique solutions \citep{bertsekas1997nonlinear}. In the case when a problem does not meet these criteria, a common strategy is to replace the original problem with a well-chosen approximation for which BCD is guaranteed to converge. A detailed discussion of such strategies can be found in \cite{razaviyayn2013unified} and \cite{yang2019inexact}. Here, we make use of the block successive upper minimization (BSUM) algorithm of \cite{razaviyayn2013unified} to ensure convergence of BCD for \eqref{eqn:KKT-EPR}, which neither has a smooth objective nor results in subproblems with unique solutions. Specifically, we make use of the result stated in Theorem \ref{thm:ProxConv}.

\begin{theorem} \label{thm:ProxConv} \citep{razaviyayn2013unified}
Let $f$ in \eqref{eqn:MCP} be nonsmooth and the subproblems \eqref{eqn:SubMCP} be convex, but not necessarily strictly convex. BCD converges to a stationary point of \eqref{eqn:MCP} if executed with the following set of subproblems (instead of \eqref{eqn:SubMCP}):
\begin{equation}
    \label{eqn:ProxSubMCP}
    \begin{aligned}
    \minimize\limits_{{x}_i \in \mathcal{X}_{i}} \quad & f(\bar{x}_1, \ldots, \bar{x}_{i-1}, x_{i}, \bar{x}_{i+1}, \ldots, \bar{x}_\ell) + \frac{1}{2 \gamma} \, \lVert x_i - \bar{x}_i \rVert_2^2 \quad \forall \, i \in \{1, \ldots, \ell\}, \\
    \end{aligned}
\end{equation}
where $\gamma$ is a scalar parameter.
\end{theorem}

Algorithm \ref{alg:BCD} shows the BCD algorithm that we use to solve \eqref{eqn:MCP}. 

\begin{algorithm}
\begin{algorithmic}[1]
    \State initialize: $x \gets \bar{x}$
    \While{convergence criteria are not satisfied}
    \ForAll{$i = 1, \dots, \ell$}
        \State solve \eqref{eqn:ProxSubMCP}, obtain $x_i^*$,  $\bar{x}_i \gets x_i^*$
    \EndFor
    \EndWhile
    \State \Return $\bar{x}$
    \caption{A cyclic BCD implementation on \eqref{eqn:MCP}.}
    \label{alg:BCD}
\end{algorithmic}
\end{algorithm}

The following corollary follows directly from Theorem \ref{thm:ProxConv}.
\begin{corollary}
Suppose $x^*$ is the accumulation point of the sequence of iterates $\bar{x}$ generated by Algorithm \ref{alg:BCD}. Then $x^*$ is a KKT point of \eqref{eqn:MCP}.
\end{corollary}

\paragraph{Initialization.}
Since our proposed solution method only guarantees to return a KKT point of \eqref{eqn:KKT-Convex}, it is important to seed it with a good initial guess. This is, in general, difficult to achieve without any knowledge of the distribution of noise in the training data. Here we detail a heuristic initialization strategy which, as evident from our extensive computational results in Section \ref{sec:CompRes}, is quite effective in finding an estimate for parameters $\theta$ and $\omega$ that perform well in predicting optimal decisions for unseen $u$ values. 

In our approach, we first initialize $\hat{x}_i$ with the noisy decision data $x_i$ for all $i \in \mathcal{I}$, i.e., $\hat{x}_0 = x$. Next, we solve the following two problems to sequentially initialize the $\hat{\omega}$ and $(\hat{\theta}, \lambda, \mu)$ parameters, respectively: 

\begin{equation}
\label{eqn:Init_Omega}
    \begin{aligned}
    \hat{\omega}_0 \in \argmin\limits_{\hat{\omega} \in \Omega} \quad &
    \sum\limits_{i \in \mathcal{I}} \mathrm{max} \{ 0, \hat{g}(\hat{x}_{0i}, u_i; \hat{\omega})\} +
    \sum\limits_{i \in \mathcal{I}} \, \lvert \hat{h}(\hat{x}_{0i}, u_i; \hat{\omega}) \rvert
    \end{aligned}
\end{equation}

\begin{equation}
\label{eqn:Init_theta}
    \begin{aligned}
    (\hat{\theta}_0, {\lambda}_0, \mu_0) \in \argmin\limits_{\hat{\theta} \in \Theta, \lambda, \mu} \quad & 
    \sum\limits_{i \in \mathcal{I}} \, \lvert \nabla \hat{f} (\hat{x}_{0i}, u_i; \hat{\theta}) + \lambda_i^{\top} \nabla \hat{g} (\hat{x}_{0i}, u_i; \hat{\omega}_0) + \mu_i^{\top} \nabla \hat{h} (\hat{x}_{0i}, u_i; \hat{\omega}_0) \rvert \, + \\
     & \sum\limits_{i \in \mathcal{I}} \, \lvert \lambda_{i}^{\top} \, \hat{g}(\hat{x}_{0i}, u_i; \hat{\omega}_0) \rvert \\
    \st \quad & \lambda_i \geq 0 \quad \forall \, i \in \mathcal{I} \\
    \end{aligned}
\end{equation}

% Note that in the case when $\omega$ are already known, i.e., all constraints are known, problem \eqref{eqn:Init_theta} is exactly the IOP formulation proposed by \cite{Boyd2011} to estimate unknown objective parameters. However, as shown by \cite{Aswani2018} and \cite{Thai2018}, this formulation does not lead to good parameter estimates when the training data is noisy.  

\paragraph{Hyperparameter Setting.} The overall algorithm involves three hyperparameters: the initial penalty parameters $c_1$, the factor $\rho$ by which $c$ are increased after every iteration of Algorithm \ref{alg:Outerloop}, and the regularization parameter $\gamma$ for BCD subproblems \eqref{eqn:ProxSubMCP}. Out of these three, we find that $c_1$ and $\rho$ have a large impact on the quality of the solution determined by our approach. While choosing $c_1$ and $\rho$ to be of large magnitude helps to quickly find a solution for which $\lVert P \rVert = 0$, very often the resulting solution is of low quality. Here, assigning small values to $c_1$ and $\rho$ generally solves the problem, but slows the convergence of Algorithm \ref{alg:Outerloop}. In our computational experiments, we utilize a validation dataset to determine sufficiently large values of $c_1$ and $\rho$ for which the algorithm converges quickly yet generates good solutions.

The regularization parameter $\gamma$ does not affect the quality of the solution, but has an impact on BCD's rate of convergence. We find that not having the regularization term in \eqref{eqn:ProxSubMCP}, i.e., setting $\gamma$ to $\infty$ may make BCD cycle between non-stationary points without converging to any particular solution. On the other hand, making $\gamma$ small stabilizes the algorithm but makes convergence extremely slow as the regularization term becomes dominant. In our computational experiments, we observe that assigning a large value to $\gamma$ such as $\gamma = 10^6$ typically works well in terms of stabilizing the BCD iterations yet ensuring fast convergence to a stationary point.

%Finally, we conclude this section by highlighting the connections between our approach and some recent NLP literature. Our overall approach using Algorithm \ref{alg:Outerloop} with Algorithm \ref{alg:BCD} as the solution method for \eqref{eqn:KKT-EPR} involves penalizing difficult coupling constraints in the objective and leveraging the resulting favorable mathematical structure by using BCD. This is similar in spirit to alternating direction method of multipliers (ADMM), which has recently become very popular for large-scale nonlinear problems \citep{Boyd2011}. However, in general, ADMM convergence is guaranteed only for convex optimization problems with linearly coupled blocks of variables. Here, we show that by using the penalty function approach and updating the penalty parameters in the outer loop, we can achieve convergence even for nonconvex optimization problems with nonlinear coupling constraints between its variables.

\section{Computational Case Studies} \label{sec:CompRes}

We apply the proposed solution methods for the data-driven IO problem to three case studies. The first case study is based on the so-called water-filling problem from the field of information theory \citep{kalpana2015}. We use this example to showcase the computational advantages of Algorithms \ref{alg:Outerloop} and \ref{alg:BCD}. The next two case studies, one addressing learning the risk preferences of a decision maker and the second related to estimating parameters of a resource allocation market, demonstrate the utility of IO in systems engineering applications. All model instances presented in this section were implemented in Julia v1.6.1 using the mathematical optimization modeling environment JuMP v0.21.10 \citep{DunningHuchetteLubin2017}. We applied Gurobi v9.1.2 to solve all convex optimization problems, and all nonconvex NLPs were solved using IPOPT v0.7.0.

\subsection{The Water-Filling Problem}
We consider the following variation of the water-filling problem from the wireless communication literature that is used to determine optimal power allocation for a multidimensional communication channel:
\begin{equation}
    \label{eqn:wat_fil}
    \begin{aligned}
    \maximize\limits_{x \in \mathbb{R}^D_+} \quad & \sum\limits_{d = 1}^{D} \theta_d \, \mathrm{log} (x_d + u_d) \\
    \st \quad & \sum\limits_{d = 1}^{D} \omega_d x_d = \omega_{D+1},
    \end{aligned}
\end{equation}
where the objective function represents the total communication rate and the constraint limits the total amount of power that can be allocated to the system. Here we assume that the weighting parameters $\theta$ for the objective and $\omega$ for the constraint are unknown. The goal is to estimate these parameters by observing changes in the optimal power allocation decisions $x$ based on fluctuations in the input parameters $u$. In the context of IO, this problem has previously been considered by \cite{Aswani2018}; however, their work is limited to the estimation of objective parameters assuming $\omega$ to be known.

For each instance of the IOP, the training data are generated as follows. We first create arbitrary vectors $\theta \in \mathbb{R}^D$ and $\{\omega_d\}_{d=1}^{D}$ with $\omega_d \in \mathbb{R}^D$ by sampling their individual elements from the uniform distribution $\mathcal{U}(1.00, 1.10)$. For all instances, we set the right-hand-side constraint parameter $\omega_{D+1}$ to the dimensionality of the problem $D$. We then sample a set of input vectors $u_i \in \mathbb{R}^D$ for each $i \in \mathcal{I}$ such that $u_{id} \sim \mathcal{U}(1.00, 2.00)$ for every $d \in \{1, \ldots, D\}$. Next, keeping $\theta$ and $\omega$ the same, we solve these $|\mathcal{I}|$ instances of \eqref{eqn:wat_fil} to obtain the optimal power allocation decisions $x_i^{*}.$ Finally, we distort these true optimal solutions by adding a Gaussian noise $\gamma \sim \mathcal{N}(0, \sigma^2 \mathbb{I})$ to obtain the noisy dataset as $x_i = x_i^* + \gamma$.   
In this study, to show the computational advantage of Algorithm \ref{alg:Outerloop}, we consider large datasets of up to 1,000 samples with FOPs of varying dimensionality. We also consider varying levels of noise in the datasets by changing the value of $\sigma$. A specific case is hence represented by $D$, $\sigma$, and $|\mathcal{I}|$, and we solve ten random instances for each case to obtain reliable performance statistics for different solution methods. While setting up the IOP, we assume $\hat{\Theta} = \left\{\hat{\theta} \, | \, 10^{-4} \leq \hat{\theta}_d \leq 10\ \, \forall \, d \in \{1, \ldots, D\} \right\}$ and $\hat{\Omega} = \left\{ \hat{\omega} \, | \, \hat{\omega}_d \geq 0\ \, \forall \, d \in \{1, \ldots, D\}, \, \hat{\omega}_{D+1} = 1 \right\}.$ The sets $\hat{\Theta}$ and $\hat{\Omega}$ have been chosen so that trivial solutions, such as $\theta = 0$, are eliminated from the feasible solution space of the IOP.

\subsubsection{Convergence criterion for Algorithm \ref{alg:Outerloop}} \label{sec:ConvCrit}
In Algorithm \ref{alg:Outerloop}, we consider the algorithm to have converged when the norm of the penalty term becomes less than a certain threshold value $\epsilon$. However, in our computational experiments, we find convergence of the algorithm with BCD in the inner loop to be very slow. Figure \ref{fig:CS1_ConvCrit}a shows the evolution of $\lVert P \rVert_1$ as a function of the number of iterations with iteration 0 representing the initialization solution. As can be observed, the value of $\lVert P \rVert_1$ reaches a peak before settling into a slow decay phase, after which the desired convergence threshold takes a long time to reach. This may, at first glance, make BCD-based Algorithm \ref{alg:Outerloop} appear an impractical approach for our problem. However, as can be seen in Figure \ref{fig:CS1_ConvCrit}b, we find that the prediction error of the generated estimates $\hat{\theta}$ and $\hat{\omega}$ stabilizes when $\lVert P \rVert_1$ enters the slow decay phase. Running the algorithm after this point does not improve the prediction accuracy. Therefore, we modify our convergence criterion to be the stabilization of the prediction error instead of $\lVert P \rVert_1 \leq \epsilon.$      

Note that the slow reduction in $\lVert P \rVert_1$ is observed only when BCD is used to solve the penalty reformulation. When solving the penalty reformulation directly with IPOPT, we find that, in most cases, the algorithm converges in just one or two iterations.

\begin{figure}
\centering
\hfill
\begin{minipage}{.45\textwidth}
\begin{tikzpicture}[scale=0.7]
    \begin{axis}[xmin = 0, xmax = 30, ylabel = {$\lVert P \rVert_1$}, xlabel = Iteration]
		\addplot[opacity = 1, blue, line width = 2] table[x index = {0}, y index = 1] {Paper_CaseStudy1_Figure3.dat};
\end{axis}
\end{tikzpicture}
\caption*{(a) Norm of the penalty term}
\end{minipage}
\begin{minipage}{.45\textwidth}
\begin{tikzpicture}[scale=0.7]
    \begin{axis}[xmin = 0, xmax = 10, ymax = 2e9, ylabel = {$d_{\mathcal{V}}(x^*, \hat{x})$}, xlabel = Iteration]
		\addplot[opacity = 1, blue, line width = 2] table[x index = {0}, y index = 1] {Paper_CaseStudy1_Figure2.dat};
\end{axis}
\end{tikzpicture}
\caption*{ (b) Prediction error on the test dataset}
\end{minipage}
\caption{Example progress of $\lVert P \rVert_1$ and $d_{\mathcal{V}}(x^*, \hat{x})$ values as the number of iterations of Algorithm \ref{alg:Outerloop} increases.}
\label{fig:CS1_ConvCrit}
\end{figure}

\subsubsection{Computational Performance}
We report results comparing the computational performance of three different solution methods for \eqref{eqn:KKT-Convex}: 
\begin{itemize}
\item Method A - solve the problem directly with IPOPT, 
\item Method B - use Algorithm \ref{alg:Outerloop} with IPOPT to solve the penalty reformulation, and 
\item Method C - use BCD (Algorithm \ref{alg:BCD}) to solve the penalty reformulation in Algorithm \ref{alg:Outerloop}. 
\end{itemize}

For a fair comparison, we initialize all three methods with the same initial solution obtained using the strategy described in Section \ref{sec:Algorithm}. All model instances are solved with a time limit of 300~s using 24 cores and 16 GB of memory on the Mesabi cluster of the Minnesota Supercomputing Institute (MSI).

Table \ref{tab:Ipopt} summarizes the results for Method A. For a specific $D$, $\sigma$, and $|\mathcal{I}|$, we show the median value of the computation time required to solve the ten random instances. An ``n/a" entry in this column indicates that IPOPT was unable to solve any of the instances in the specified time limit. Additionally, we provide the median, maximum, and minimum prediction errors of the generated estimates $\hat{\theta}$ and $\hat{\omega}$. To obtain these statistics, along with every instance of training data, we also generate a test dataset $\mathcal{V}$ of 100 inputs. This test data consists of $(u, x^*)$ pairs, where $u$ values are generated in the same manner as for the training data, and $x^*$ are the corresponding true optimal solutions obtained using the same $\theta$ and $\omega$ as the ones used to generate the noisy training data. Once $\hat{\theta}$ and $\hat{\omega}$ have been found, we use them to solve the problems in the test dataset and evaluate the prediction error as the total Manhattan distance between the true solutions and their estimated values, i.e., we use the  metric $d_{\mathcal{V}}(x^*, \hat{x}) = \sum\limits_{v \in \mathcal{V}} \lVert x_v^* - \hat{x}_v \rVert_{1}$. 

\begin{table}[h!]
\resizebox{\textwidth}{!}{%
\centering
\begin{tabular}{@{}ccccccccccc@{}}
\toprule
  \multirow{3}{*}[-2em]{$\sigma$} &
  \multirow{3}{*}[-2em]{$\mathcal{I}$} &
  \multicolumn{4}{c}{$D = 50$} &
  \multicolumn{4}{c}{$D = 100$} \\ \cmidrule(lr){3-6} \cmidrule(lr){7-10} 
   &
   &
  \multirow{2}{2.4cm}[-0.2em] {\centering Median Computation Time (s)} &
  \multicolumn{3}{c}{Prediction Error} &
  \multirow{2}{2.4cm}[-0.2em]{\centering Median Computation Time (s)} &
  \multicolumn{3}{c}{Prediction Error}  \\         \addlinespace[3pt]
 \cmidrule(lr){4-6} \cmidrule(lr){8-10}
&     &    & median & min & max &    & median & min & max &      \\  \midrule
\multirow{2}{*}{0.01} & $D$   & 36 & 62.96           & 26.90       & 1,702.86         & 96 & 580.16           & 79.21       & 1,305.69   \\
& 1,000 & 198 & 5.61          & 4.70      & 11.33       & 206 & 6.26           & 5.56       & 6.54    \\
\multirow{2}{*}{0.05} & $D$  & 64 & 1.97         & 0.29      & 16.62       & 63 & 3.10 $\times 10^5$           & 118.57      & 2.93 $\times 10^8$       \\
& 1,000 & n/a &     -     &    -    &   -     & n/a &     -      &  -     &   -    \\
\multirow{2}{*}{0.1}  & $D$  & 60  & 2.33            & 0.54      & 7.07       & 64 & 6.34 $\times 10^6$          & 6.13 $\times 10^4$      & 8.11 $\times 10^7$   \\
& 1,000 & n/a  &   -        &   -    &    -    & n/a &    -       & -      & - \\ \bottomrule
 \end{tabular}}
\caption{\label{tab:Ipopt} Computational performances of IPOPT (Method A) on instances of \eqref{eqn:KKT-Convex} designed to estimate parameters of \eqref{eqn:wat_fil}.}
\end{table}

From the data in Table \ref{tab:Ipopt}, one can observe that IPOPT can solve smaller instances very quickly but struggles with larger problems. Moreover, even for problems with smaller datasets, there is a large variability in the prediction accuracy values obtained with estimated parameters. The solved instances with $D = 50$ generally exhibit good prediction performance, in contrast to instances with $D = 100$ irrespective of the dataset size. We find that, generally, IPOPT is not a robust solution method for difficult instances of \eqref{eqn:KKT-Convex} with either high dimensionality or large datasets.

Next, using the same metrics as defined for Method A, we show a comparison of the performances of solution methods B and C in Table \ref{tab:BSUM}. While implementing Algorithm \ref{alg:Outerloop} on IOP instances of \eqref{eqn:wat_fil}, we set $c_1 = 500 e$, where $e$ is a vector of all ones and $\rho = 1,000$. These values were obtained by trial-and-error as described in Section \ref{sec:Algorithm}.

\begin{table}[h!]
\resizebox{\textwidth}{!}{%
\centering
\begin{tabular}{@{}ccccccccccc@{}}
\toprule
\multirow{3}{*}[-3em]{$D$} &
  \multirow{3}{*}[-3em]{$\sigma$} &
  \multirow{3}{*}[-3em]{$\mathcal{I}$} &
  \multicolumn{4}{c}{Solved with IPOPT (Method B)} &
  \multicolumn{4}{c}{Solved with BCD (Method C)} \\ \cmidrule(lr){4-7} \cmidrule(lr){8-11} 
 &
   &
   &
  \multirow{2}{2.4cm}[-1em] {\centering Median Computation Time (s)} &
  \multicolumn{3}{c}{Prediction Error} &
  \multirow{2}{2.4cm}[-1em]{\centering Median Computation Time* (s)} &
  \multicolumn{3}{c}{Prediction Error}  \\         \addlinespace[3pt]
 \cmidrule(lr){5-7} \cmidrule(lr){9-11}
                     &                       &     &    & median & min & max &    & median & min & max      \\ \addlinespace[10pt] \midrule
\multirow{6}{*}{50} & \multirow{2}{*}{0.01} & 50   & 11 & 64.20           & 26.90       & 104.08         & 1 & 41.60           & 23.65       & 141.88   \\
                     &                       & 1,000 & 219 & 1.85          & 1.70      & 18.49       & 33 & 2.96           & 2.32       & 3.14    \\
                     & \multirow{2}{*}{0.05} & 50  & 24 & 108.21         & 94.65      & 121.49       & 1 & 116.55           & 89.62      & 255.47      \\
                     &                       & 1,000 & 300  & 71.83         & 9.25       & 6,533.90       & 31 & 14.77          & 10.99      & 15.65        \\
                     & \multirow{2}{*}{0.1}  & 50  & 55  & 174.99            & 143.54      & 286.04       & 2 & 166.19          & 134.18      & $1.02 \times 10^7$   \\
                     &                       & 1,000 & 300  & 4,504.87          & 3,994.46      & 4,512.91       & 61 & 37.81          & 31.18      & 44.19       \\ \midrule
\multirow{6}{*}{100} & \multirow{2}{*}{0.01} & 100   & 101  & 102.24          & 42.67       & 797.99        & 4 & 109.16          & 42.21     & 338.88   \\
                     &                       & 1,000 & 300 & 22.60          & 19.73      & 160.78       & 51 & 5.88          & 5.71      & 6.53       \\
                     & \multirow{2}{*}{0.05} & 100  & 113  & 155.36         & 135.21      & 1,426.01       & 7  & 156.68         & 134.54      & $1.44 \times 10^6$      \\
                     &                       & 1,000 & 300  & $2.32 \times 10^5$         &  $1.36 \times 10^5$     & $6.13 \times 10^5 $        & 112 & 47.85           & 27.95      & 53.28    \\
                     & \multirow{2}{*}{0.1}  & 100  & 109  & 195.23 & 91.05 & 265.68                        & 10 & 248.60         & 199.15      & $2.18 \times 10^8$       \\
                     &                       & 1,000 & 300  & $4.4 \times 10^5$         & $8.6 \times 10^4$      & $6.43 \times 10^5$       & 107  & 70.60            & 65.86      & 93.36       \\ \bottomrule
\end{tabular}}
\caption{\label{tab:BSUM} Comparison of computational performances of solution methods B and C on an IOP based on random instances of \eqref{eqn:wat_fil}. *Based on the modified convergence criterion described in Section \ref{sec:ConvCrit}.}
\end{table}

The computational advantage of using a BCD-based decomposition strategy is immediately apparent from the median computation times in Table \ref{tab:BSUM}. Even for the smaller instances where we do not expect decomposition to vastly outperform a full-space formulation, BCD outputs a solution as much as ten times faster compared to when IPOPT is used to solve the penalty reformulation of \eqref{eqn:KKT-Convex}. For larger problems, we observe that while IPOPT starts reaching the time limit of 300 s, BCD is able to output a solution in approximately 100 s. 

\subsubsection{Solution Quality}

Next, we focus on the quality of the solutions generated by Methods B and C. For smaller instances, we find that both methods yield solutions of similar prediction accuracy. However, as we move to problems with larger datasets, IPOPT does not produce a good solution in any of the problem instances. In contrast, with BCD, we see that the median prediction error stays low for all cases. The only case where BCD \emph{may} result in a bad solution is when the training data has a high level of noise yet a small dataset is provided to generate an estimate of the missing FOP parameters. 
%We suspect this is because the approximate nature of BCD makes it more prone to the presence of local minima whose number is higher when the dataset is smaller. 
Finally, Figure \ref{fig:CS1_StatCons} shows that the solutions obtained from BCD asymptotically converge towards the one with minimum prediction error. This empirically confirms that our algorithm generates statistically consistent solutions.

\usepgfplotslibrary{colorbrewer}
% initialize Set1-4 from colorbrewer (we're comparing 4 classes),
\pgfplotsset{compat = 1.15, cycle list/Set1-8} 
% Tikz is loaded automatically by pgfplots
\usetikzlibrary{pgfplots.statistics, pgfplots.colorbrewer} 

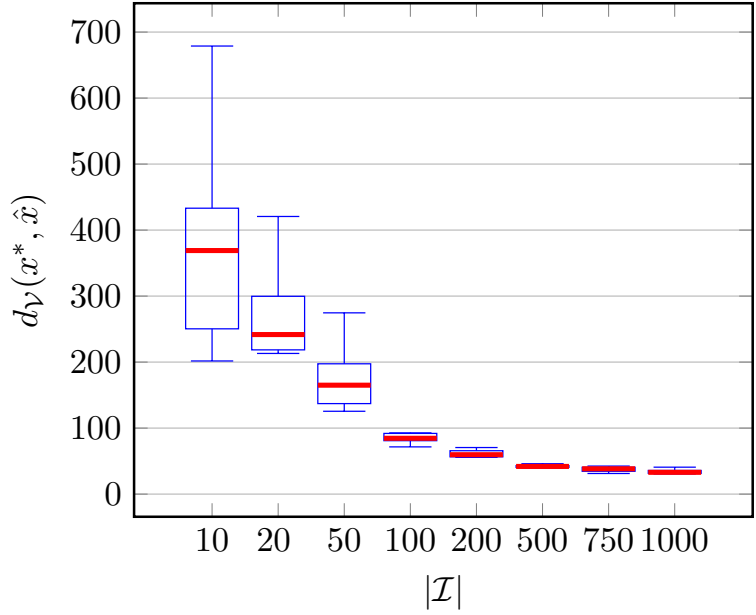
\begin{figure}[h!]
\centering
    \begin{tikzpicture}[scale=1.2]
    	\pgfplotstableread[col sep=comma]{CaseStudy1_Figure1.csv}{\csvdata}
    	% Boxplot groups columns, but we want rows
    	\pgfplotstabletranspose\datatransposed{\csvdata} 
    	\begin{axis}[
    	   % cycle list/Dark2,
    		boxplot/draw direction = y,
    		x axis line style = {opacity=100, thick},
    		y axis line style = {thick},
    		% axis x line* = bottom,
    		% axis y line = left,
    		enlarge y limits,
    		ymajorgrids,
    		xtick = {1, 2, 3, 4, 5, 6, 7, 8},
    		xticklabel style = {align=center},
    		xticklabels = {$10$,$20$ , $50$, $100$, $200$, $500$, $750$, $1000$},
    		xtick style = {}, % Hide tick line
    		ylabel = {$d_{\mathcal{V}}(x^*, \hat{x})$},
    		ytick = {0, 100, 200, 300, 400, 500, 600, 700},
    		xlabel = {$|\mathcal{I}|$}
    	]
    		\foreach \n in {1,...,8} {
    			\addplot+[boxplot = {every median/.style={red,ultra thick}}, draw=blue] table[y index=\n] {\datatransposed};
    		}
    	\end{axis}
    \end{tikzpicture}
    \caption{Change in prediction error as a function of the training dataset size. The box plot shows the interquartile ranges of prediction error values for the ten random instances solved with $D = 100$ and $\sigma = 0.1$ using solution method C.}
    \label{fig:CS1_StatCons}
\end{figure}

Overall, we find that Algorithm \ref{alg:Outerloop} is a significantly more robust solution method for problem \eqref{eqn:KKT-Convex} compared to a direct application of standard commercial NLP solvers. Moreover, when the penalty reformulation of \eqref{eqn:KKT-Convex} is an MCP (as in this case), using BCD in the inner loop of the algorithm produces marked improvement in its computational efficiency without compromising (or rather enhancing) the final solution quality.

\subsection{Estimating Risk Preferences}

When people make decisions under significant uncertainty, their choices strongly depend on their individual risk preferences. Knowing these risk preferences can be very helpful to higher-level decision makers; for example, companies may be able to create better products and services for their risk-averse customers, and policymakers could develop more effective regulations and incentives for their communities. However, often these risk preferences are not explicitly known but need to be estimated from observed decisions.

In risk-averse optimization, the decision-making problem is commonly formulated as a two-stage stochastic program of the following form:
\begin{equation}
\label{eqn:RiskProblem}
    \minimize_{x \in \mathcal{X}} \; f_1(x) + \mathbb{E}(Q(x,\xi)) + \lambda \, \rho(Q(x,\xi)),
\end{equation}
where $x$ and $\xi$ denote the first-stage decisions and uncertain parameters (or random variables), respectively, $\mathcal{X}$ is the first-stage feasible set, and $f_1(x)$ is the first-stage cost function. The recourse function is denoted by $Q(x,\xi)$ and is defined as
\begin{equation}
    Q(x,\xi) = \min_{y \in \mathcal{Y}(x,\xi)} \; f_2(x,y,\xi)
\end{equation}
with $\mathcal{Y}(x,\xi)$ being the set of feasible second-stage decisions $y$ and $f_2(x,y,\xi)$ being the second-stage cost function. The risk function is denoted by $\rho$, and the objective is to minimize a weighted sum of the expected cost and risk, where the weighting factor $\lambda$ is chosen by the decision maker based on their level of risk aversion.

Many common risk functions are convex, which allows us to learn them using the proposed data-driven IO framework, given that the overall FOP is also convex. Indeed, in the risk literature, there is the notion of \textit{coherent} risk measures, which exhibit properties that are deemed desirable for the purpose of quantifying risk \citep{Artzner1999}; one of those properties is convexity. In the following example, we consider a simplified version of the risk-based electricity procurement problem introduced by \citet{Zhang2016}, where the conditional value-at-risk (CVaR), one of the most popular coherent risk measures, is used as the risk function $\rho$.

In this problem, the decision maker needs to schedule the operation of a production plant that consumes a large amount of electricity, which can be purchased from a power contract or from the spot market. The prices for the power contract are known, but one must commit to the purchase amount in advance, e.g. one day or one week ahead. Spot electricity prices, on the other hand, are uncertain and only become known shortly before the time of delivery. While purchasing electricity from the power contract bears no risk, the average contract price is typically higher than the average spot price, and it leaves less room for taking advantage of fluctuations in the spot price; hence, there is a trade-off between purchasing from the power contract and the spot market, which one seeks to optimize. Assuming that the uncertain spot prices follow a discrete probability distribution, the problem can be formulated as the following deterministic equivalent of the corresponding two-stage stochastic program:
\begin{subequations}
\label{eqn:risk_pref}
\begin{align}
    \minimize_{x,y,q,v,w} \quad & \sum_{t \in \mathcal{T}} c_t \, x_t^2 + \sum_{s \in \mathcal{S}} p_s \sum_{t \in \mathcal{T}} \left( b_t \, q_{ts} + r_{ts} \, y_{ts} \right) + \lambda \left( v + \frac{1}{1-\alpha} \sum_{s \in \mathcal{S}} p_s \, w_s \right) \label{eqn:ElecProc_obj} \\
    \st \quad & i^{\min} \leq i^0 + \sum_{t'=1}^t \left( q_{ts} - d_t \right) \leq i^{\max} \quad \forall \, t \in \mathcal{T}, \, s \in \mathcal{S} \label{eqn:ElecProc_inventory} \\
    & \sum_{t \in \mathcal{T}} \left( q_{ts} - d_t \right) \geq 0 \quad \forall \, s \in \mathcal{S} \label{eqn:ElecProc_totalP} \\
    & m \, q_{ts} = x_t + y_{ts} \quad \forall \, t \in \mathcal{T}, \, s \in \mathcal{S} \label{eqn:ElecProc_consumption} \\
    & \sum_{t \in \mathcal{T}} \left( b_t \, q_{ts} + r_{ts} \, y_{ts} \right) - v \leq w_s \quad \forall \, s \in \mathcal{S} \label{eqn:ElecProc_CVaR} \\
    & x_t \in \mathbb{R}_+ \quad \forall \, t \in \mathcal{T} \\
    & q_{ts}, \, y_{ts} \in \mathbb{R}_+ \quad \forall \, t \in \mathcal{T}, \, s \in \mathcal{S} \\
    & w_s \in \mathbb{R}_+ \quad \forall \, s \in \mathcal{S},
\end{align}
\end{subequations}
where $\mathcal{T}$ is the set of time periods that defines the scheduling horizon, and $\mathcal{S}$ is the set of possible uncertainty realizations (or scenarios). Each scenario $s$ is defined by its probability $p_s$ and spot price profile $r_{ts}$. The first-stage decisions are the amounts of electricity purchased from the power contract, denoted by $x_t$, across the scheduling horizon. The second-stage variables are the production amounts, $q_{ts}$, and the amounts of electricity purchased from the spot market, $y_{ts}$. The cost of electricity purchased from the power contract in time period $t$ is $c_t \, x_t^2$, and the production cost (excluding the electricity cost) is assumed to be a linear function of the production amount defined by the coefficients $b_t$. In constraints \eqref{eqn:ElecProc_inventory}, the parameters $i^0$, $i^{\min}$, $i^{\max}$, and $d_t$ denote the initial product inventory level, the minimum inventory, the maximum inventory, and the product demand in time period $t$, respectively. Constraints \eqref{eqn:ElecProc_totalP} ensure that the total amount of product produced over the course of the scheduling horizon is not less than the total demand. We assume that the electricity consumption is a linear function of the production amount defined by the coefficient $m$ and, per constraints $\eqref{eqn:ElecProc_consumption}$, must be equal to the total amount of electricity procured in each time period. The term in the objective function \eqref{eqn:ElecProc_obj} that is multiplied by the weighting factor $\lambda$, along with constraints \eqref{eqn:ElecProc_CVaR}, represents the $\alpha$-CVaR. For a given $\alpha \in (0,1)$, the $\alpha$-CVaR is defined as the expected cost greater than the $\alpha$-VaR, which is the $\alpha$-quantile of the cost distribution. Note the use of the auxiliary variables $v$ and $w_s$, where for each scenario in which the cost is greater than $v$, $w_s$ takes the value of the difference between the cost and $v$.

In this case study, we assume that the decision maker must decide on their electricity procurement strategy two days in advance. We discretize the 48-hour time horizon into one-hour time periods, yielding $\lvert \mathcal{T} \rvert = 48$. We assume that the decision maker considers 35 different spot electricity price scenarios and $\alpha = 0.9$ when making their decisions. We set up the IOP by assuming that we have access to historical electricity procurement data in which the decision maker made decisions under changing sets of spot electricity price scenarios. Our goal here is to use our knowledge of price scenarios and the corresponding procurement decisions to estimate the decision maker's risk tolerance ($\lambda$) and preferred safety stock level ($i^{\mathrm{min}}$), where the latter can also be seen as an indicator of the decision maker's level of risk aversion. All other parameters in \eqref{eqn:risk_pref} are assumed to be known. We assume that only the optimal values of $(x_t, q_{ts})$ for all $t \in \mathcal{T}$ and \emph{one of the $s \in \mathcal{S}$} can be observed  since in practice, we are only able to see the decisions made in the scenario that actually realized. As such, this case study is an example of an IOP in which only a subset of the decision variables of the FOP can be observed.

We provide the values of all the parameters for \eqref{eqn:risk_pref} including the contract and spot electricity price profiles in the supplementary material. To generate each training data point in $\mathcal{I}$, we begin by selecting 35 spot price scenarios at random from the 40 included in the supplementary material. We then solve \eqref{eqn:risk_pref} with these scenarios for certain $\lambda$ and $i^{\mathrm{min}}$ to obtain the optimal decisions. We select one of the $s$ in $\mathcal{S}$ at random as the scenario that actually realizes, and we store the corresponding optimal $x_t$ and $q_{st}$ values for all $t \in \mathcal{T}$, as well as the $r_{st}$ values for all $t \in \mathcal{T}$ and $s \in \mathcal{S}$, as a data point in $\mathcal{I}$. 

Here, the FOP is a QP with 3,444 variables and 8,553 constraints, making the IOP a large-scale nonconvex NLP. As a result, we use BCD to solve the penalty reformulation in Algorithm \ref{alg:Outerloop}. This algorithm was implemented with $c_1 = 0.01e$ and $\rho = 0.01$. Because we only observe a subset of the decision variables, the initialization strategy discussed in Section \ref{sec:Algorithm} cannot be applied here, so the algorithm is initialized with random values for the missing parameters. We solve the IOP using a training dataset of 5 partial observations (i.e., $\lvert \mathcal{I} \rvert = 5$), which we find is sufficient to accurately estimate the values of $\lambda$ and $i^{\mathrm{min}}$. Furthermore, we solve five instances of the IOP, each with a different randomly generated training dataset, to confirm that our results are robust against price scenario selection. 

The results of our computational experiments are summarized in Table \ref{tab:cs2_results}. We consider three different cases involving decision makers with varying risk tolerances. In Case 1, where the decision maker is almost risk-neutral (with a low $\lambda$ value), our algorithm predicts the same $\lambda$ and $i^{\mathrm{min}}$ values as the ones used to generate the five instances of training data. However, in the other two cases we observe some discrepancies between predicted and true values. In Case 2, where the decision maker has a medium risk tolerance, the predicted value of $i^{\mathrm{min}}$ is greater than its true value. This is because the constraints and parameters of \eqref{eqn:risk_pref} are such that the product inventory can never fall below 50 and therefore, our algorithm ``thinks" that 50 is the safety stock level for this decision maker. In Case 3, where the large $\lambda$ indicates a highly risk-averse decision maker, the reason for the difference in the $\lambda$ values is that for any $\lambda \gtrsim 5$, the decision maker is so risk-averse that it always chooses to buy all the electricity from the power contract. It is important to note that in Cases 2 and 3, despite the differences in the true and predicted parameter values, solving \eqref{eqn:risk_pref} with any of the predicted $\lambda$ and $i^{\mathrm{min}}$ values will always yield the same decisions as the ones obtained using the true FOP. In other words, in all three cases, our solution approach provides parameters that effectively minimize the decision error on both seen and unseen spot price scenarios which is what the objective of \eqref{eqn:IOP} is.

\begin{table}[h!]
\centering
\begin{tabular}{@{}cccccccc@{}}
\toprule
\multirow{3}{*}{} & \multicolumn{1}{c}{\multirow{3}{*}{}} & \multicolumn{1}{c}{\multirow{3}{*}[-1em]{True Values}} & \multicolumn{5}{c}{ Predicted Values} \\ \cmidrule(lr){4-8} 
                  & \multicolumn{1}{c}{} & \multicolumn{1}{c}{} & \multicolumn{5}{c}{Instances}         \\  
                  & \multicolumn{1}{c}{} & \multicolumn{1}{c}{} & 1     & 2     & 3     & 4     & 5     \\ \midrule
\multirow{2}{*}{Case 1} &        $\lambda$              & 0.1                  & 0.099 & 0.1   & 0.099 & 0.1   & 0.099 \\
                  &        $i^{\mathrm{min}}$              & 75                   & 74.99 & 74.99 & 74.99 & 74.99 & 74.99 \\ \midrule
\multirow{2}{*}{Case 2} &      $\lambda$                & 2                    & 1.97  & 1.98  & 1.97  & 1.99  & 1.965 \\
                  &       $i^{\mathrm{min}}$               & 40                   & 50    & 50    & 50    & 50    & 50    \\ \midrule
\multirow{2}{*}{Case 3} &      $\lambda$                & 7                    & 6.98  & 19.52 & 5.52  & 5.36  & 23.28 \\
                  &        $i^{\mathrm{min}}$              & 60                   & 59.99 & 59.99 & 60    & 59.99 & 60    \\ \bottomrule
\end{tabular}
\caption{\label{tab:cs2_results} Results of computational experiments based on electricity procurement problem \eqref{eqn:risk_pref}.}
\end{table}

In this problem, we only had access to a subset of the decisions from the model, yet we find that even with partial observability of the decisions, we were able to correctly recover the unknown parameters. However, this is not always the case. In most situations, partial observability of the decisions will expand the solution space of the inverse problem which may lead to estimated parameter values that do not predict the decisions accurately. This highlights the need to develop conditions under which partially observable decisions can be used to accurately recover the parameters of the FOP with some theoretical guarantee.

\subsection{Resource Allocation Market}
In this case study, we consider a market with a set of buyers $\mathcal{B}$ that enter a bargaining game to acquire a limited set of resources $\mathcal{G}$. Each buyer has a utility $u_b$ as a function of resource allocations and a disagreement point $d_b$, which is the status quo that player $b$ will revert to if no agreement is reached. We assume that the buyers agree to cooperate with each other such that the resource allocation is \emph{fair}. \citet{nash1950bargaining} showed that if the utility set is compact and convex, solving the following problem achieves the unique solution that satisfies certain axioms, which represent a notion of fairness:
\begin{equation}
\label{eqn:nbg}
\begin{aligned}
    \maximize_x \quad & \prod_{b \in \mathcal{B}} (u_b(x) - d_b) \\
    \st \quad & p(x) \leq 0,
\end{aligned}
\end{equation}
where $x$ are the resource allocation decisions and $p$ is a convex function encoding the constraints on the market. Problem \eqref{eqn:nbg} is not convex; however, a logarithmic transformation of the objective function allows it to be reformulated as a convex optimization problem. 

Here, we assume that the market is constrained by two sets of constraints: first, buyers are limited in the amount of resources they can acquire, and second, the total amount of each resource available in the market is limited. A convexified reformulation of \eqref{eqn:nbg} for this market is as follows:
\begin{equation}
\tag{RAMP}
\label{eqn:res_allocation}
\begin{aligned}
\maximize\limits_{x \in \mathbb{R}^{|\mathcal{B}| \times |\mathcal{G}|}_{+}} \quad & \sum\limits_{b \in \mathcal{B}} \mathrm{log} \left( \sum\limits_{g \in \mathcal{G}} x_{bg} \right) \\
\st \quad & \sum\limits_{g \in \mathcal{G}} p_{bg} x_{bg}^2 \leq m_b \quad \forall \, b \in \mathcal{B} \\
& \sum\limits_{b \in \mathcal{B}} x_{bg} \leq c_g \quad \forall \, g \in \mathcal{G}, \\
\end{aligned}
\end{equation}
where $x_{bg}$ refers to the amount of a resource $g$ that a buyer $b$ ends up acquiring from the market. The first set of constraints place limits on how much resource a buyer can acquire and the second set of constraints model the capacity of the market. The objective function in this model seeks to determine a fair allocation of utilities, $u_b = \sum\limits_{g \in \mathcal{G}}x_{bg}$ for each buyer $b \in \mathcal{B}$ assuming that at the disagreement point, $d_b = 0$.

For our IO set up, we consider a scenario where the resource providers cooperate among themselves to learn the constraints of the buyers, i.e., the resource providers are interested in estimating the values of parameters $p_{bg}$ and $m_b$ for each $b \in \mathcal{B}.$ To learn these parameters, the resource providers rely on their knowledge of historical capacity fluctuation data (changes in the values of $c_g$) and the corresponding fair allocations. We assume that this historical resource allocation data is only available with high noise.

To generate training data for this case study, we assume a market with 20 buyers and 5 resources. For each instance of the IOP, we sample individual elements of $p \in \mathbb{R}^{|\mathcal{B}| \times |\mathcal{G}|}$ and $m \in \mathbb{R}^{\lvert \mathcal{B} \rvert}$ randomly from uniform distributions $\mathcal{U}(0.5, 2)$ and $\mathcal{U}(0.5, 1)$, respectively. We then sample a set of capacity vectors $c_{ig} \in \mathbb{R}^{\lvert \mathcal{G} \rvert}$ for each $i \in \mathcal{I}$ such that each of their elements follow the distribution $\mathcal{U}(0.1, 1.5)$. Next, keeping $p$ and $m$ the same, we solve the $\lvert \mathcal{I} \rvert$ instances of \eqref{eqn:res_allocation} corresponding to each of the $\lvert \mathcal{I} \rvert$ $c_{ig}$ vectors. We then distort the obtained optimal resource allocation decisions $x^*_{ibg}$ by adding a Gaussian noise $\gamma \sim \mathcal{N}(0, \sigma^2 \mathbb{I})$ to them. Here we find that with our chosen parameter values, the mean value of $x^*_{ibg}$ for each $i \in \mathcal{I}$ results in being approximately $0.3$, and therefore to simulate a high noise scenario, we set the value of $\sigma$ to $0.05$.

Note that the single-level IOP reformulation \eqref{eqn:KKT-Convex} for \eqref{eqn:res_allocation} is not an MCP. Therefore, we are restricted to Algorithm \ref{alg:Outerloop} to solve this problem. Nonetheless, we find that the constraint parameters of \eqref{eqn:res_allocation} can be learned with just a few data points, making the problem size amenable for Algorithm \ref{alg:Outerloop}. While learning the parameters, we assume that the $m_b$ values for all $b \in \mathcal{B}$ are 1 and instead learn the normalized $p_{bg}$ parameters, i.e., the $p_{bg}/m_b$ values. In Figure \ref{fig:CS3_param_fit}, we compare the accuracy of estimates generated through IO with simple regression-based estimates. To generate the regression estimates, we apply the least-squares method on the noisy data to find the best-fit parameters for the quadratic constraints assuming that they hold as equalities in all cases. While we agree that assuming that these constraints always hold as equalities is a rather strong assumption and may negatively affect the quality of the estimates, this is typically how this problem will be solved if no other information is available. Indeed, as we show in Figure \ref{fig:CS3_param_fit}, IO can generate more accurate estimates with just 5 data points than what regression can produce with 100 data points. On the one hand, this is an obvious result because our assumption of a fair allocation provides additional information about the data. On the other hand, it also shows the power of the proposed framework in being highly data-efficient by allowing the incorporation of domain knowledge into the inverse problem.

\usepgfplotslibrary{colorbrewer}
\pgfplotsset{compat = 1.15, cycle list/Set1-9} 
% Tikz is loaded automatically by pgfplots
\usetikzlibrary{pgfplots.statistics, pgfplots.colorbrewer} 
\pgfplotsset{tick label style={font=\large},
    label style={font=\large}}

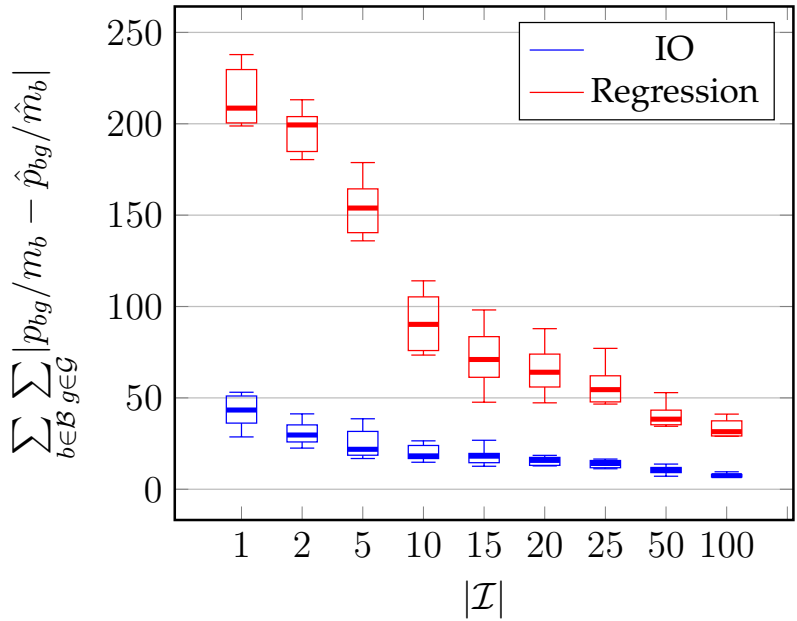
\begin{figure}[ht]
\centering
    \begin{tikzpicture}[scale=1.2]
    	\pgfplotstableread[col sep=comma]{IO_error.csv}{\csvdataio}     
         \pgfplotstableread[col sep=comma]{regr_error.csv}{\csvdataregr}     
    	% Boxplot groups columns, but we want rows
    \pgfplotstabletranspose\datatransposedio{\csvdataio} 
    \pgfplotstabletranspose\datatransposedregr{\csvdataregr} 
    	\begin{axis}[
    	   % cycle list/Dark2,
            name=border,
    		boxplot/draw direction = y,
    		x axis line style = {opacity=100, thick},
    		y axis line style = {thick},
    		% axis x line* = bottom,
    		% axis y line = left,
    		enlarge y limits,
    		ymajorgrids,
    		xtick = {0.5, 1.5, 2.5, 3.5, 4.5, 5.5, 6.5, 7.5, 8.5, 9.5},
    		xticklabel style = {align=center},
    		xticklabels = {$1$,$2$,$5$,$10$,$15$,$20$,$25$, $50$,$100$},
    		xtick style = {}, % Hide tick line
    		ylabel = {$\sum\limits_{b \in \mathcal{B}} \sum\limits_{g \in \mathcal{G}} \lvert p_{bg}/m_{b} - \hat{p}_{bg}/\hat{m}_b \rvert$},
    		xlabel = {$|\mathcal{I}|$},
            legend entries={IO,Regression},
            legend pos=north east]
    	]
        \addlegendimage{no markers, blue}
        \addlegendimage{no markers, red}
    	\foreach \n in {1,...,9} {
            \addplot+[boxplot = {every median/.style={blue,ultra thick}, draw position=0.5 + \n - 1, box extend = 0.5}, draw=blue] table[y index=\n] {\datatransposedio};
            % \ifthenelse{\n=1}{\addlegendentry{IO}}{}
            \addplot+[boxplot = {every median/.style={red,ultra thick}, draw position=0.5 + \n - 1, box extend = 0.5}, draw=red] table[y index=\n] {\datatransposedregr};
            % \ifthenelse{\n=1}{\addlegendentry{Regression}}{}
    	}
    	\end{axis}
    \end{tikzpicture}
    \caption{Change in the error in the normalized parameter estimate values as a function of the training dataset size. Box plots show the interquartile ranges of error values for ten random instances of the problem generated by following the scheme described in the text.}
    \label{fig:CS3_param_fit}
\end{figure}

\section{Conclusions} \label{sec:Conclusions}
In this work, we described the idea of discovering unknown decision-making processes from observed decisions using mathematical optimization as a natural model for decision making. Using this concept, inferring the decision-making model is equivalent to estimating the unknown parameters of the underlying optimization problem, which is referred to as data-driven inverse optimization. We considered the case in which the underlying decision-making process can be formulated as a convex optimization problem. We formulated the IO problem as a bilevel program with as many lower-level problems as the number of available observations of prior decisions. To address the computational challenge associated with solving large instances of the IO problem, we proposed an efficient penalty-based BCD algorithm that leverages the decomposable structure of the problem.

We conducted extensive computational experiments to benchmark the performance of the proposed solution method against standard commercial solvers. In large instances, we show that using our BCD algorithm does not only reduce the computation time but also results in higher-quality solutions. Furthermore, we present two additional computational case studies based on practically relevant problems, one concerned with estimating risk preferences and the other one aimed at uncovering local constraint parameters in a multiplayer Nash bargaining formulation. Here, our results indicate that IO can recover highly accurate estimates of the parameters of interest while using only a very small number of data points.

Finally, we would like to highlight that there are many important directions to consider for future work including decomposition algorithms for the problems where BCD cannot be applied, dealing with partially observable decisions, incorporating prior beliefs, and using adaptive sampling to reduce the amount of data required to learn the parameters.

\section*{Acknowledgments}
The authors gratefully acknowledge the financial support from the National Science Foundation under Grant \#2044077 as well as the Minnesota Supercomputing Institute (MSI) at the University of Minnesota for providing resources that contributed to the research results reported in this paper. R.G. acknowledges financial support from a departmental fellowship sponsored by 3M and a Doctoral Dissertation Fellowship from the University of Minnesota. 

\bibliographystyle{newapa}
% Compile offline later with \bibliographystyle{newapa}, don't know why Overleaf can't find that style file.
\bibliography{library}

\begin{thebibliography}{}

\bibitem[\protect\citeauthoryear{Abbeel \& Ng}{Abbeel \& Ng}{2004}]{Abbeel2004}
Abbeel, P. \& Ng, A.~Y. (2004).
\newblock {Apprenticeship learning via inverse reinforcement learning}.
\newblock In {\em Proceedings of the 21st International Conference on Machine
  Learning}.

\bibitem[\protect\citeauthoryear{Ahuja \& Orlin}{Ahuja \&
  Orlin}{2001}]{Ahuja2001}
Ahuja, R.~K. \& Orlin, J.~B. (2001).
\newblock {Inverse Optimization}.
\newblock {\em Operations Research}, {\em 49\/}(5), 771--783.

\bibitem[\protect\citeauthoryear{Akhtar, Kolarijani \& Esfahani}{Akhtar
  et~al.}{2022}]{Akhtar2022}
Akhtar, S.~A., Kolarijani, A.~S., \& Esfahani, P.~M. (2022).
\newblock {Learning for Control: An Inverse Optimization Approach}.
\newblock {\em IEEE Control Systems Letters}, {\em 6\/}(c), 187--192.

\bibitem[\protect\citeauthoryear{Anitescu}{Anitescu}{2000}]{anitescu2000nonlinear}
Anitescu, M. (2000).
\newblock Nonlinear programs with unbounded lagrange multiplier sets.
\newblock {\em Preprint ANL/MCS-P796-0200, Mathematics and Computer Science
  Division, Argonne National Laboratory, Argonne, IL}.

\bibitem[\protect\citeauthoryear{Anitescu}{Anitescu}{2005}]{anitescu2005using}
Anitescu, M. (2005).
\newblock On using the elastic mode in nonlinear programming approaches to
  mathematical programs with complementarity constraints.
\newblock {\em SIAM Journal on Optimization}, {\em 15\/}(4), 1203--1236.

\bibitem[\protect\citeauthoryear{Artzner, Delbaen, Eber \& Heath}{Artzner
  et~al.}{1999}]{Artzner1999}
Artzner, P., Delbaen, F., Eber, J.~M., \& Heath, D. (1999).
\newblock {Coherent measures of risk}.
\newblock {\em Mathematical Finance}, {\em 9\/}(3), 203--228.

\bibitem[\protect\citeauthoryear{Aswani, Shen \& Siddiq}{Aswani
  et~al.}{2018}]{Aswani2018}
Aswani, A., Shen, Z.-J.~M., \& Siddiq, A. (2018).
\newblock {Inverse Optimization with Noisy Data}.
\newblock {\em Operations Research}, {\em 66\/}(3), 870--892.

\bibitem[\protect\citeauthoryear{Babier, Chan, Lee, Mahmood \& Terekhov}{Babier
  et~al.}{2021}]{Babier2021}
Babier, A., Chan, T. C.~Y., Lee, T., Mahmood, R., \& Terekhov, D. (2021).
\newblock {An Ensemble Learning Framework for Model Fitting and Evaluation in
  Inverse Linear Optimization}.
\newblock {\em INFORMS Journal on Optimization}, {\em 3\/}(2), 119--138.

\bibitem[\protect\citeauthoryear{Bal{\'{a}}zsi, {Van Oudenaarden} \&
  Collins}{Bal{\'{a}}zsi et~al.}{2011}]{Balazsi2011}
Bal{\'{a}}zsi, G., {Van Oudenaarden}, A., \& Collins, J.~J. (2011).
\newblock {Cellular decision making and biological noise: From microbes to
  mammals}.
\newblock {\em Cell}, {\em 144\/}(6), 910--925.

\bibitem[\protect\citeauthoryear{Ba{\~{n}}ares-Alc{\'{a}}ntara, Westerberg \&
  Rychener}{Ba{\~{n}}ares-Alc{\'{a}}ntara et~al.}{1985}]{Banares-Alcantara1985}
Ba{\~{n}}ares-Alc{\'{a}}ntara, R., Westerberg, A.~W., \& Rychener, M.~D.
  (1985).
\newblock {Development of an expert system for physical property predictions}.
\newblock {\em Computers and Chemical Engineering}, {\em 9\/}(2), 127--142.

\bibitem[\protect\citeauthoryear{Bangi \& Kwon}{Bangi \&
  Kwon}{2020}]{Bangi2020}
Bangi, M. S.~F. \& Kwon, J. S.~I. (2020).
\newblock {Deep hybrid modeling of chemical process: Application to hydraulic
  fracturing}.
\newblock {\em Computers and Chemical Engineering}, {\em 134}.

\bibitem[\protect\citeauthoryear{Benson, Sen \& Shanno}{Benson
  et~al.}{2007}]{benson2009convergence}
Benson, H.~Y., Sen, A., \& Shanno, D.~F. (2007).
\newblock Convergence analysis of an interior-point method for nonconvex
  nonlinear programming.
\newblock {\em Available on Optimization Online}.

\bibitem[\protect\citeauthoryear{Benson, Sen, Shanno \& Vanderbei}{Benson
  et~al.}{2006}]{benson2006interior}
Benson, H.~Y., Sen, A., Shanno, D.~F., \& Vanderbei, R.~J. (2006).
\newblock Interior-point algorithms, penalty methods and equilibrium problems.
\newblock {\em Computational Optimization and Applications}, {\em 34\/}(2),
  155--182.

\bibitem[\protect\citeauthoryear{Bertsekas}{Bertsekas}{1997}]{bertsekas1997nonlinear}
Bertsekas, D.~P. (1997).
\newblock Nonlinear programming.
\newblock {\em Journal of the Operational Research Society}, {\em 48\/}(3),
  334--334.

\bibitem[\protect\citeauthoryear{Bertsimas, Gupta \& Paschalidis}{Bertsimas
  et~al.}{2015}]{Bertsimas2015}
Bertsimas, D., Gupta, V., \& Paschalidis, I.~C. (2015).
\newblock {Data-driven estimation in equilibrium using inverse optimization}.
\newblock {\em Mathematical Programming}, {\em 153\/}(2), 595--633.

\bibitem[\protect\citeauthoryear{Birge, Hortacsu \& Pavlin}{Birge
  et~al.}{2017}]{Birge2017}
Birge, J.~R., Hortacsu, A., \& Pavlin, M. (2017).
\newblock {Inverse Optimization for the Recovery of Market Structure from
  Market Outcomes: An Application to the MISO Electricity Market}.
\newblock {\em Operations Research}, {\em 65\/}(4), 837--855.

\bibitem[\protect\citeauthoryear{Bollas, Barton \& Mitsos}{Bollas
  et~al.}{2009}]{Bollas2009}
Bollas, G.~M., Barton, P.~I., \& Mitsos, A. (2009).
\newblock {Bilevel optimization formulation for parameter estimation in
  vapor-liquid(-liquid) phase equilibrium problems}.
\newblock {\em Chemical Engineering Science}, {\em 64\/}(8), 1768--1783.

\bibitem[\protect\citeauthoryear{Boukouvala \& Floudas}{Boukouvala \&
  Floudas}{2017}]{Boukouvala2017}
Boukouvala, F. \& Floudas, C.~A. (2017).
\newblock {ARGONAUT: AlgoRithms for Global Optimization of coNstrAined grey-box
  compUTational problems}.
\newblock {\em Optimization Letters}, {\em 11\/}(5), 895--913.

\bibitem[\protect\citeauthoryear{Boyd, Parikh, Chu, Peleato \& Eckstein}{Boyd
  et~al.}{2011}]{Boyd2011}
Boyd, S., Parikh, N., Chu, E., Peleato, B., \& Eckstein, J. (2011).
\newblock {Distributed Optimization and Statistical Learning via the
  Alternating Direction Method of Multipliers}.
\newblock {\em Foundations and Trends in Machine Learning}, {\em 3\/}(1),
  1--122.

\bibitem[\protect\citeauthoryear{Burgard \& Maranas}{Burgard \&
  Maranas}{2003}]{Burgard2003}
Burgard, A.~P. \& Maranas, C.~D. (2003).
\newblock {Optimization-based framework for inferring and testing hypothesized
  metabolic objective functions}.
\newblock {\em Biotechnology and Bioengineering}, {\em 82\/}(6), 670--677.

\bibitem[\protect\citeauthoryear{Burton \& Toint}{Burton \&
  Toint}{1992}]{Burton1992}
Burton, D. \& Toint, P.~L. (1992).
\newblock {On an instance of the inverse shortest paths problem}.
\newblock {\em Mathematical Programming}, {\em 53\/}(1-3), 45--61.

\bibitem[\protect\citeauthoryear{Chan \& Kaw}{Chan \& Kaw}{2019}]{Chan2019a}
Chan, T.~C. \& Kaw, N. (2019).
\newblock {Inverse optimization for the recovery of constraint parameters}.
\newblock {\em European Journal of Operational Research}, {\em 282\/}(2),
  415--427.

\bibitem[\protect\citeauthoryear{Chan, Lee \& Terekhov}{Chan
  et~al.}{2019}]{Chan2019}
Chan, T.~C., Lee, T., \& Terekhov, D. (2019).
\newblock {Inverse optimization: Closed-form solutions, geometry, and goodness
  of fit}.
\newblock {\em Management Science}, {\em 65\/}(3), 1115--1135.

\bibitem[\protect\citeauthoryear{Chan, Craig, Lee \& Sharpe}{Chan
  et~al.}{2014}]{Chan2014}
Chan, T. C.~Y., Craig, T., Lee, T., \& Sharpe, M.~B. (2014).
\newblock {Generalized Inverse Multiobjective Optimization with Application to
  Cancer Therapy}.
\newblock {\em Operations Research}, {\em 62\/}(3), 680--695.

\bibitem[\protect\citeauthoryear{Chan, Eberg, Forster, Holloway, Ieraci,
  Shalaby \& Yousefi}{Chan et~al.}{2021}]{Chan2021a}
Chan, T. C.~Y., Eberg, M., Forster, K., Holloway, C., Ieraci, L., Shalaby, Y.,
  \& Yousefi, N. (2021).
\newblock {An Inverse Optimization Approach to Measuring Clinical Pathway
  Concordance}.
\newblock {\em Management Science}, (January 2022).

\bibitem[\protect\citeauthoryear{Chan, Mahmood \& Zhu}{Chan
  et~al.}{2021}]{Chan2021}
Chan, T. C.~Y., Mahmood, R., \& Zhu, I.~Y. (2021).
\newblock
\newblock {Inverse Optimization: Theory and Applications}, 1--71.

\bibitem[\protect\citeauthoryear{Choromanska, Henaff, Mathieu, Arous \&
  LeCun}{Choromanska et~al.}{2015}]{choromanska2015loss}
Choromanska, A., Henaff, M., Mathieu, M., Arous, G.~B., \& LeCun, Y. (2015).
\newblock The loss surfaces of multilayer networks.
\newblock In {\em Artificial intelligence and statistics}, (pp.\ 192--204).
  PMLR.

\bibitem[\protect\citeauthoryear{Chow, Ritchie \& Jeong}{Chow
  et~al.}{2014}]{Chow2014}
Chow, J.~Y., Ritchie, S.~G., \& Jeong, K. (2014).
\newblock {Nonlinear inverse optimization for parameter estimation of
  commodity-vehicle-decoupled freight assignment}.
\newblock {\em Transportation Research Part E: Logistics and Transportation
  Review}, {\em 67}, 71--91.

\bibitem[\protect\citeauthoryear{Danilova, Dvurechensky, Gasnikov, Gorbunov,
  Guminov, Kamzolov \& Shibaev}{Danilova et~al.}{2022}]{danilova2022}
Danilova, M., Dvurechensky, P., Gasnikov, A., Gorbunov, E., Guminov, S.,
  Kamzolov, D., \& Shibaev, I. (2022).
\newblock Recent theoretical advances in non-convex optimization.
\newblock In {\em High-Dimensional Optimization and Probability}  (pp.\
  79--163). Springer.

\bibitem[\protect\citeauthoryear{Dempe, Kalashnikov, P{\'e}rez-Vald{\'e}s \&
  Kalashnykova}{Dempe et~al.}{2015}]{dempe2015bilevel}
Dempe, S., Kalashnikov, V., P{\'e}rez-Vald{\'e}s, G.~A., \& Kalashnykova, N.
  (2015).
\newblock Bilevel programming problems.
\newblock {\em Energy Systems. Springer, Berlin}.

\bibitem[\protect\citeauthoryear{Dunning, Huchette \& Lubin}{Dunning
  et~al.}{2017}]{DunningHuchetteLubin2017}
Dunning, I., Huchette, J., \& Lubin, M. (2017).
\newblock Jump: A modeling language for mathematical optimization.
\newblock {\em SIAM Review}, {\em 59\/}(2), 295--320.

\bibitem[\protect\citeauthoryear{Gautam \& Seider}{Gautam \&
  Seider}{1979}]{Gautam1979}
Gautam, R. \& Seider, W.~D. (1979).
\newblock {Computation of phase and chemical equilibrium: Part I. Local and
  constrained minima in Gibbs free energy}.
\newblock {\em AIChE Journal}, {\em 25\/}(6), 991--999.

\bibitem[\protect\citeauthoryear{Ghobadi \& Mahmoudzadeh}{Ghobadi \&
  Mahmoudzadeh}{2020}]{Ghobadi2020}
Ghobadi, K. \& Mahmoudzadeh, H. (2020).
\newblock {Multi-point Inverse Optimization of Constraint Parameters}.
\newblock {\em arXiv preprint}, {\em arXiv:2001}.

\bibitem[\protect\citeauthoryear{Glass, Aigner, Viell, Jupke \& Mitsos}{Glass
  et~al.}{2017}]{Glass2017a}
Glass, M., Aigner, M., Viell, J., Jupke, A., \& Mitsos, A. (2017).
\newblock {Liquid-liquid equilibrium of 2-methyltetrahydrofuran/water over wide
  temperature range: Measurements and rigorous regression}.
\newblock {\em Fluid Phase Equilibria}, {\em 433}, 212--225.

\bibitem[\protect\citeauthoryear{Glass \& Mitsos}{Glass \&
  Mitsos}{2018}]{Glass2018}
Glass, M. \& Mitsos, A. (2018).
\newblock {Parameter estimation in reactive systems subject to sufficient
  criteria for thermodynamic stability}.
\newblock {\em Chemical Engineering Science}, {\em 197}, 420--431.

\bibitem[\protect\citeauthoryear{Gupta \& Zhang}{Gupta \&
  Zhang}{2022}]{gupta2022deco}
Gupta, R. \& Zhang, Q. (2022).
\newblock Decomposition and adaptive sampling for data-driven inverse linear
  optimization.
\newblock {\em INFORMS Journal on Computing}.

\bibitem[\protect\citeauthoryear{Heuberger}{Heuberger}{2004}]{Heuberger2004}
Heuberger, C. (2004).
\newblock {Inverse combinatorial optimization: A survey on problems, methods,
  and results}.
\newblock {\em Journal of Combinatorial Optimization}, {\em 8\/}(3), 329--361.

\bibitem[\protect\citeauthoryear{Hey \& Orme}{Hey \& Orme}{1994}]{Hey1994}
Hey, J.~D. \& Orme, C. (1994).
\newblock {Investigating generalizations of expected utility theory using
  experimental data}.
\newblock {\em Econometrica}, {\em 62\/}(6), 1291--1326.

\bibitem[\protect\citeauthoryear{Iraj \& Terekhov}{Iraj \&
  Terekhov}{2021}]{Iraj2021}
Iraj, E.~H. \& Terekhov, D. (2021).
\newblock
\newblock {Comparing Inverse Optimization and Machine Learning Methods for
  Imputing a Convex Objective Function}.

\bibitem[\protect\citeauthoryear{Iyengar \& Kang}{Iyengar \&
  Kang}{2005}]{Iyengar2005}
Iyengar, G. \& Kang, W. (2005).
\newblock {Inverse conic programming with applications}.
\newblock {\em Operations Research Letters}, {\em 33\/}(3), 319--330.

\bibitem[\protect\citeauthoryear{Jain \& Kar}{Jain \& Kar}{2017}]{jain2017non}
Jain, P. \& Kar, P. (2017).
\newblock Non-convex optimization for machine learning.
\newblock {\em arXiv preprint arXiv:1712.07897}.

\bibitem[\protect\citeauthoryear{Jin, Netrapalli, Ge, Kakade \& Jordan}{Jin
  et~al.}{2021}]{Jin2021}
Jin, C., Netrapalli, P., Ge, R., Kakade, S.~M., \& Jordan, M.~I. (2021).
\newblock On nonconvex optimization for machine learning: Gradients,
  stochasticity, and saddle points.
\newblock {\em Journal of the ACM (JACM)}, {\em 68\/}(2), 1--29.

\bibitem[\protect\citeauthoryear{Kalpana \& Khan}{Kalpana \&
  Khan}{2015}]{kalpana2015}
Kalpana, N. \& Khan, M. Z.~A. (2015).
\newblock Fast computation of generalized waterfilling problems.
\newblock {\em IEEE Signal Processing Letters}, {\em 22\/}(11), 1884--1887.

\bibitem[\protect\citeauthoryear{Keshavarz, Wang \& Boyd}{Keshavarz
  et~al.}{2011}]{Keshavarz2011}
Keshavarz, A., Wang, Y., \& Boyd, S. (2011).
\newblock {Imputing a Convex Objective Function}.
\newblock In {\em IEEE International Symposium on Intelligent Control}, (pp.\
  613--619).

\bibitem[\protect\citeauthoryear{McCarl, Moskowitz \& Furtan}{McCarl
  et~al.}{1977}]{mccarl1977}
McCarl, B.~A., Moskowitz, H., \& Furtan, H. (1977).
\newblock Quadratic programming applications.
\newblock {\em Omega}, {\em 5\/}(1), 43--55.

\bibitem[\protect\citeauthoryear{McFarland}{McFarland}{1977}]{McFarland1977}
McFarland, D.~J. (1977).
\newblock {Decision making in animals}.
\newblock {\em Nature}, {\em 269\/}(5623), 15--21.

\bibitem[\protect\citeauthoryear{{Mohajerin Esfahani}, Shafieezadeh-Abadeh,
  Hanasusanto \& Kuhn}{{Mohajerin Esfahani} et~al.}{2018a}]{Esfahani2018}
{Mohajerin Esfahani}, P., Shafieezadeh-Abadeh, S., Hanasusanto, G.~A., \& Kuhn,
  D. (2018a).
\newblock {Data-driven inverse optimization with imperfect information}.
\newblock {\em Mathematical Programming}, {\em 167\/}(1), 191--234.

\bibitem[\protect\citeauthoryear{{Mohajerin Esfahani}, Shafieezadeh-Abadeh,
  Hanasusanto \& Kuhn}{{Mohajerin Esfahani}
  et~al.}{2018b}]{MohajerinEsfahani2018}
{Mohajerin Esfahani}, P., Shafieezadeh-Abadeh, S., Hanasusanto, G.~A., \& Kuhn,
  D. (2018b).
\newblock {Data-driven inverse optimization with imperfect information}.
\newblock {\em Mathematical Programming}, {\em 167\/}(1), 191--234.

\bibitem[\protect\citeauthoryear{Morgenstern \& von Neumann}{Morgenstern \& von
  Neumann}{1944}]{Morgenstern1944}
Morgenstern, O. \& von Neumann, J. (1944).
\newblock {\em {Theory of Games and Economic Behavior}}.
\newblock Princeton University Press.

\bibitem[\protect\citeauthoryear{Nash~Jr}{Nash~Jr}{1950}]{nash1950bargaining}
Nash~Jr, J.~F. (1950).
\newblock The bargaining problem.
\newblock {\em Econometrica: Journal of the econometric society}, 155--162.

\bibitem[\protect\citeauthoryear{Nocedal \& Wright}{Nocedal \&
  Wright}{2006}]{nocedal2006numerical}
Nocedal, J. \& Wright, S. (2006).
\newblock {\em Numerical optimization}.
\newblock Springer Science \& Business Media.

\bibitem[\protect\citeauthoryear{Parker \& Smith}{Parker \&
  Smith}{1990}]{Parker1990}
Parker, G.~A. \& Smith, J.~M. (1990).
\newblock {Optimality theory in evolutionary biology}.
\newblock {\em Nature}, {\em 348\/}(6296), 27--33.

\bibitem[\protect\citeauthoryear{Razaviyayn, Hong \& Luo}{Razaviyayn
  et~al.}{2013}]{razaviyayn2013unified}
Razaviyayn, M., Hong, M., \& Luo, Z.-Q. (2013).
\newblock A unified convergence analysis of block successive minimization
  methods for nonsmooth optimization.
\newblock {\em SIAM Journal on Optimization}, {\em 23\/}(2), 1126--1153.

\bibitem[\protect\citeauthoryear{Rich \& Venkatasubramanian}{Rich \&
  Venkatasubramanian}{1987}]{Rich1987}
Rich, S.~H. \& Venkatasubramanian, V. (1987).
\newblock {Model-based reasoning in diagnostic expert systems for chemical
  process plants}.
\newblock {\em Computers and Chemical Engineering}, {\em 11\/}(2), 111--122.

\bibitem[\protect\citeauthoryear{Rich \& Venkatasubramanian}{Rich \&
  Venkatasubramanian}{1989}]{Rich1989}
Rich, S.~H. \& Venkatasubramanian, V. (1989).
\newblock {Causality‐based failure‐driven learning in diagnostic expert
  systems}.
\newblock {\em AIChE Journal}, {\em 35\/}(6), 943--950.

\bibitem[\protect\citeauthoryear{R{\"{o}}nnqvist, Svenson, Flisberg \&
  J{\"{o}}nsson}{R{\"{o}}nnqvist et~al.}{2017}]{Ronnqvist2017}
R{\"{o}}nnqvist, M., Svenson, G., Flisberg, P., \& J{\"{o}}nsson, L.~E. (2017).
\newblock {Calibrated route finder: Improving the safety, environmental
  consciousness, and cost effectiveness of truck routing in Sweden}.
\newblock {\em Interfaces}, {\em 47\/}(5), 372--395.

\bibitem[\protect\citeauthoryear{Rosen}{Rosen}{1967}]{Rosen1967}
Rosen, R. (1967).
\newblock {\em {Optimality Principles in Biology}}.
\newblock Springer.

\bibitem[\protect\citeauthoryear{Rossi, Cardozo-Filho \& Guirardello}{Rossi
  et~al.}{2009}]{Rossi2009}
Rossi, C.~C., Cardozo-Filho, L., \& Guirardello, R. (2009).
\newblock {Gibbs free energy minimization for the calculation of chemical and
  phase equilibrium using linear programming}.
\newblock {\em Fluid Phase Equilibria}, {\em 278\/}(1-2), 117--128.

\bibitem[\protect\citeauthoryear{Saez-Gallego, Morales, Zugno \&
  Madsen}{Saez-Gallego et~al.}{2016}]{Saez-Gallego2016}
Saez-Gallego, J., Morales, J.~M., Zugno, M., \& Madsen, H. (2016).
\newblock {A Data-Driven Bidding Model for a Cluster of Price-Responsive
  Consumers of Electricity}.
\newblock {\em IEEE Transactions on Power Systems}, {\em 31\/}(6), 5001--5011.

\bibitem[\protect\citeauthoryear{Sammut, Hurst, Kedzier \& Michie}{Sammut
  et~al.}{1992}]{Sammut1992}
Sammut, C., Hurst, S., Kedzier, D., \& Michie, D. (1992).
\newblock {Learning to fly}.
\newblock {\em Machine Learning Proceedings}, 385--393.

\bibitem[\protect\citeauthoryear{Scheel \& Scholtes}{Scheel \&
  Scholtes}{2000}]{scheel2000mathematical}
Scheel, H. \& Scholtes, S. (2000).
\newblock Mathematical programs with complementarity constraints: Stationarity,
  optimality, and sensitivity.
\newblock {\em Mathematics of Operations Research}, {\em 25\/}(1), 1--22.

\bibitem[\protect\citeauthoryear{Schoemaker}{Schoemaker}{1991}]{Schoemaker1991}
Schoemaker, P. J.~H. (1991).
\newblock {The quest for optimality: A positive heuristic of science?}
\newblock {\em Behavioral and Brain Sciences}, {\em 14}, 205--245.

\bibitem[\protect\citeauthoryear{Shahmoradi \& Lee}{Shahmoradi \&
  Lee}{2021}]{Shahmoradi2021a}
Shahmoradi, Z. \& Lee, T. (2021).
\newblock {Quantile Inverse Optimization: Improving Stability in Inverse Linear
  Programming}.
\newblock {\em Operations Research}, (January 2022).

\bibitem[\protect\citeauthoryear{Shen, Diamond, Udell, Gu \& Boyd}{Shen
  et~al.}{2017}]{shen2017disciplined}
Shen, X., Diamond, S., Udell, M., Gu, Y., \& Boyd, S. (2017).
\newblock Disciplined multi-convex programming.
\newblock In {\em 2017 29th Chinese Control And Decision Conference (CCDC)},
  (pp.\ 895--900). IEEE.

\bibitem[\protect\citeauthoryear{Steels}{Steels}{1995}]{Steels1995}
Steels, L. (1995).
\newblock {When are robots intelligent autonomous agents?}
\newblock {\em Robotics and Autonomous Systems}, {\em 15\/}(1-2), 3--9.

\bibitem[\protect\citeauthoryear{Stephanopoulos}{Stephanopoulos}{1990}]{Stephanopoulos1990}
Stephanopoulos, G. (1990).
\newblock {Artificial intelligence in process engineering-current state and
  future trends}.
\newblock {\em Computers and Chemical Engineering}, {\em 14\/}(11), 1259--1270.

\bibitem[\protect\citeauthoryear{Uygun, Matthew \& Huang}{Uygun
  et~al.}{2007}]{uygun2007}
Uygun, K., Matthew, H.~W., \& Huang, Y. (2007).
\newblock Investigation of metabolic objectives in cultured hepatocytes.
\newblock {\em Biotechnology and bioengineering}, {\em 97\/}(3), 622--637.

\bibitem[\protect\citeauthoryear{W{\"a}chter \& Biegler}{W{\"a}chter \&
  Biegler}{2006}]{wachter2006implementation}
W{\"a}chter, A. \& Biegler, L.~T. (2006).
\newblock On the implementation of an interior-point filter line-search
  algorithm for large-scale nonlinear programming.
\newblock {\em Mathematical programming}, {\em 106\/}(1), 25--57.

\bibitem[\protect\citeauthoryear{Wilson \& Sahinidis}{Wilson \&
  Sahinidis}{2019}]{Wilson2019}
Wilson, Z.~T. \& Sahinidis, N.~V. (2019).
\newblock {Automated learning of chemical reaction networks}.
\newblock {\em Computers and Chemical Engineering}, {\em 127}, 88--98.

\bibitem[\protect\citeauthoryear{Wright}{Wright}{2015}]{Wright2015}
Wright, S.~J. (2015).
\newblock Coordinate descent algorithms.
\newblock {\em Mathematical Programming}, {\em 151\/}(1), 3--34.

\bibitem[\protect\citeauthoryear{Yang, Pesavento, Luo \& Ottersten}{Yang
  et~al.}{2019}]{yang2019inexact}
Yang, Y., Pesavento, M., Luo, Z.-Q., \& Ottersten, B. (2019).
\newblock Inexact block coordinate descent algorithms for nonsmooth nonconvex
  optimization.
\newblock {\em IEEE Transactions on Signal Processing}, {\em 68}, 947--961.

\bibitem[\protect\citeauthoryear{Zhang \& Liu}{Zhang \& Liu}{1996}]{Zhang1996}
Zhang, J. \& Liu, Z. (1996).
\newblock {Calculating some inverse linear programming problems}.
\newblock {\em Journal of Computational and Applied Mathematics}, {\em
  72\/}(2), 261--273.

\bibitem[\protect\citeauthoryear{Zhang, Liu \& Ma}{Zhang
  et~al.}{1996}]{Zhang1996b}
Zhang, J., Liu, Z., \& Ma, Z. (1996).
\newblock {On the inverse problem of minimum spanning tree with partition
  constraints}.
\newblock {\em Mathematical Methods of Operations Research}, {\em 44\/}(2),
  171--187.

\bibitem[\protect\citeauthoryear{Zhang, Cremer, Grossmann, Sundaramoorthy \&
  Pinto}{Zhang et~al.}{2016}]{Zhang2016}
Zhang, Q., Cremer, J.~L., Grossmann, I.~E., Sundaramoorthy, A., \& Pinto, J.~M.
  (2016).
\newblock {Risk-based integrated production scheduling and electricity
  procurement for continuous power-intensive processes}.
\newblock {\em Computers and Chemical Engineering}, {\em 86}, 90--105.

\bibitem[\protect\citeauthoryear{Zhao, Stettner, Reznik, Segre \&
  Paschalidis}{Zhao et~al.}{2015}]{Zhao2015}
Zhao, Q., Stettner, A., Reznik, E., Segre, D., \& Paschalidis, I.~C. (2015).
\newblock {Learning cellular objectives from fluxes by inverse optimization}.
\newblock {\em Proceedings of the IEEE Conference on Decision and Control},
  {\em 54rd IEEE\/}(Cdc), 1271--1276.

\end{thebibliography}

\end{document}